\documentclass[11pt]{article} 

\usepackage[utf8]{inputenc} 

\usepackage{geometry} 
\usepackage{fullpage}

\usepackage{esvect, upgreek}

\usepackage{graphicx} 

\usepackage[parfill]{parskip} 

\usepackage{tikz}

\usepackage[colorlinks=true, pdfstartview=FitV, linkcolor=blue, 
            citecolor=blue, urlcolor=blue]{hyperref}

\usepackage{booktabs} 
\usepackage{array} 
\usepackage{paralist} 
\usepackage{verbatim} 
\usepackage{subfig} 

\usepackage{fancyhdr} 
\usepackage{sectsty}
\allsectionsfont{\sffamily\mdseries\upshape} 

\usepackage[nottoc,notlof,notlot]{tocbibind} 
\usepackage[titles,subfigure]{tocloft} 

\usepackage{authblk}

\usepackage{amsmath,amsthm,amssymb}
\numberwithin{equation}{subsection}

\newtheorem{thm}{Theorem}[subsection]
\newtheorem{prop}[thm]{Proposition}

\newtheorem{lemma}[thm]{Lemma}
\newtheorem{rk}[thm]{Remark}

\theoremstyle{definition}
\newtheorem{definition}[thm]{Definition}
\newtheorem{exa}[thm]{Example}

\def\respt {\ensuremath \mathtt{r} }

\def\frb{\ensuremath{\mathfrak{b}}}
\def\dyn{\ensuremath{\mathrm{Dyn}}}

\def\fdeg{\ensuremath {\mathrm{fdeg}}}

\def\Res{\ensuremath {\mathrm{Res}}}

\def\tisfI{\ensuremath {\tilde{\mathsf{I}} }  }

\def\dem{\ensuremath {\delta_-}}
\def\dep{\ensuremath {\delta_+}}
\def\epm{\ensuremath {\epsilon_-}}
\def\epp{\ensuremath {\epsilon_+}}
\def\kam{\ensuremath {\kappa_-}}
\def\kap{\ensuremath {\kappa_+}}
\def\emm{\ensuremath {m_-}}
\def\emp{\ensuremath {m_+}}
\def\demc{\ensuremath {\delta_-^c}}
\def\depc{\ensuremath {\delta_+^c}}
\def\epmc{\ensuremath {\epsilon_-^c}}
\def\eppc{\ensuremath {\epsilon_+^c}}

\def\ZZ{\ensuremath {\mathbf{Z}}}
\def\QQ{\ensuremath {\mathbf{Q}}}
\def\RR{\ensuremath {\mathbf{R}}}
\def\CC{\ensuremath {\mathbf{C}}}

\def\FF{\ensuremath {\mathbf{F}}}
\def\KK{\ensuremath {k_{nr}} }

\def\aaa{\ensuremath {\mathrm{A}}}
\def\bbb{\ensuremath {\mathrm{B}}}
\def\ccc{\ensuremath {\mathrm{C}}}
\def\ddd{\ensuremath {\mathrm{D}}}
\def\eee{\ensuremath {\mathrm{E} }}

\def\weylh{\ensuremath {\mathcal{W}}}
\def\tts{\ensuremath {\mathtt{S}}}
\def\zm{\ensuremath {\mathtt{m}}}
\def\ha{\ensuremath {\mathcal{H}}}
\def\IM{\ensuremath {\mathrm{IM}}}
\def\calr{\ensuremath {\mathcal{R}}}

\def\PH{\ensuremath {\mathcal{P}}}

\def\utype{\ensuremath {\mathfrak{s}}}

\def\Frob{\ensuremath {\mathrm{Frob}}}

\def\End{\ensuremath {\mathrm{End}}}

\newcommand{\cind}{\textup{Ind}}

\def\parax{\ensuremath {\mathcal{X}}}
\def\parav{\ensuremath {\mathcal{V}}}

\DeclareMathOperator{\rmone}{\mathrm{I}}
\DeclareMathOperator{\rmtwo}{\mathrm{II}}
\DeclareMathOperator{\rmthree}{\mathrm{III}}
\DeclareMathOperator{\rmfour}{\mathrm{IV}}
\DeclareMathOperator{\rmfive}{\mathrm{V}}
\DeclareMathOperator{\rmsix}{\mathrm{VI}}

\title{A Note on The Spectral Transfer Morphism for Affine Hecke Algebras}

\author{Yongqi Feng \thanks{yongqifeng@science.ru.nl}}
\affil{Institute for Mathematics, Astrophysics and Particle Physics\\Radboud Universiteit Nijmegen, the Netherlands}

\date{} 

\begin{document}
\maketitle

\begin{abstract}
Lusztig \cite{Lusztig} defined and classified unipotent representations for adjoint simple groups over a non-archimedean field which are inner forms of split groups. His classification appears as a ``geometric/arithmetic correspondence'' and fits into Langlands's philosophy. Opdam \cite{Opds} introduced the notion of spectral transfer morphisms of affine Hecke algebras to study the formal degree of a unipotent discrete series representation. Based on the uniqueness property of supercuspidl unipotent representations established in \cite{FeOp}, Opdam \cite{Opdl} proved that the unipotent discrete series representations of classical groups, can be classified by the associated formal degrees, in the same spirit as Reeder's result \cite{Reexceptional} for split exceptional adjoint groups. 

The present paper aims at verifying that three finite morphisms of algebraic tori are spectral transfer morphisms, completing the proof of the main theorem in \cite{Opdl}.
\end{abstract}

\tableofcontents

\section{Introduction}

Let $k$ be a non-archimedean local field with a finite residue field $\mathbf{F}_q$. Fix a separable closure $k_s$ of $k$. Throughout this paper, assume that $\mathbf{G}$ is a connected absolutely simple algebraic group of adjoint type defined and quasi-split over $k$ such that $\mathbf{G}$ splits over $\KK$. Lusztig \cite{Lusztig} defined unipotent representations for $\mathbf{G}(k)$ (the group of $k$-rational points of $\mathbf{G}$). We recall this definition in Section \ref{sec:uptrep}. In particular, the irreducible summands of the representation $V:=\cind_{\PH_k}^{\mathbf{G}(k)} E$ of $\mathbf{G}(k)$ which is compactly induced from a cuspidal unipotent representation $(\sigma, E)$ of a parahoric subgroup $\PH_k \subset \mathbf{G}(k)$, are supercuspidal unipotent representations, provided that $\PH_k$ is a maximal parahoric subgroup of $\mathbf{G}(k)$. 

A conjecture of Hiraga, Ichino and Ikeda \cite{HII} states that for a discrete series representation of a connected reductive group defined over a local field, the formal degree $\fdeg(\pi)$ equals to $|\gamma(\varphi_\pi)|$ times a rational factor, where $\gamma(\varphi_\pi)$ is the adjoint gamma factor associated with the Langlands parameter $\varphi_\pi$ corresponding to $\pi$. An earlier result of Reeder \cite{Reexceptional} on unipotent discrete series representations of split exceptional adjoint groups is compatible with this conjecture (see \cite[\S 3.4]{HII}). 

Recently, Opdam \cite[Thm.~4.11]{Opdl} proved the Hiraga--Ichino--Ikeda (HII) conjecture for unipotent discrete series representations of $\mathbf{G}(k)$, for classical groups $\mathbf{G}$ satisfying our assumption. The main tool Opdam used is called \emph{Spectral Transfer Morphisms} for affine Hecke algebras which he developed in \cite{Opds}.


It is known that there is a bijection between the collection $\mathrm{Irr}_{upt}(\mathbf{G}(k), \PH_k, \sigma)$ of representations of $\mathbf{G}(k)$ which are irreducible summands of the compactly induced representation $V$ and simple $\ha$--modules, where $\ha:=\mathrm{End}_{\mathbf{G}(k)} (V)$. This endomorphism algebra is an extension of an affine Hecke algebra $\ha(\mathbf{G}(k), \PH_k, \sigma)$ whose parameters are explicitly determined by $\mathbf{G}$, the parahoric subgroup $\PH_k$ and the cuspidal unipotent character $\sigma$. (An explicit expression of these parameters can be found in \cite[\S 5.12]{Lusztig}.) 

If the unipotent discrete series representation is supercuspidal, the spectral transfer morphism is particularly simple. We illustrate it briefly. Suppose now $\PH_k$ is a \emph{maximal} parahoric subgroup and $\sigma$ a cuspidal unipotent character of the reductive quotient $\overline{\PH_k}$ of $\PH_k$ (we lift it to a representation of $\PH$ via the natural projection $\PH_k \twoheadrightarrow \overline{\PH_k}$). If $\pi$ belongs to $\mathrm{Irr}_{upt}(\mathbf{G}(k), \PH_k, \sigma)$, then $\pi$ is supercuspidal and $\ha(\mathbf{G}(k), \PH_k, \sigma)$ is simply $\CC$. In this case, the spectral transfer morphism has domain $\ha(\mathbf{G}(k), \PH_k, \sigma)=\CC$ and target $\ha^{\IM}$ (the Iwahori--Hecke algebra $\ha^{\IM}$ attached to $\mathbf{G}(k)$), and is given by a morphism of the associated algebraic tori 
$$
\mathcal{S}: \mathtt{Spec}\, \CC \to \mathtt{Spec}\, Z(\ha^{\IM}),
$$ 
where $Z(\ha^{\IM})$ is the centre of $\ha^{\IM}$. This morphism is determined once its image is specified. To this end, let 
$$
\varphi_\pi: \Frob^{\ZZ} \times \mathbf{SL}_2(\CC) \to {}^L G = G^\vee \rtimes \Frob^\ZZ
$$ 
be the unramified Langlands parameter associated with $\pi$. Observe that $\varphi_\pi$ is determined by the images 
\begin{equation}\label{eq:imagephi}
s:=\varphi_\pi(\Frob, I_2), \;  c:= \varphi_\pi (\mathrm{id}, \mathrm{diag}(v, v^{-1})) \text{ and }\;  u:= \varphi_\pi(\mathrm{id}, \left( \begin{smallmatrix} 1 & 1 \\ 0 & 1 \end{smallmatrix}\right) ).
\end{equation} 
It is known that if $\varphi_\pi$ is the Langlands parameter of a supercuspidal unipotent representation, then $u$ must be a distinguished unipotent element in the connected centraliser $M:=Z_{G^\vee}(s)^\circ$. The condition that $u \in M$ is distinguished unipotent forces that $s_0$ must be an isolated element in ${}^L G$ (see \cite[\S 3.7]{Retorsion}). 

Recall Borel's theorem \cite[Prop.~6.7]{Borel}, which can be formulated as: the  $G^\vee$-conjugacy classes of semi-simple elements in ${}^L G$ is bijectively corresponding to the orbit of the complex torus $\mathtt{Spec} Z(\ha^{\IM})$ under the associated finite Weyl group $W_0$. We require that the image of $\mathcal{S}$ is mapped to the $\mathrm{Ad}(G^\vee)$-orbit of $s$ in ${}^L G$ by Borel's bijection. This condition defines $\mathcal{S}$ uniquely. 

We will discuss the supercuspidal unipotent representations of non-split but quasi-split groups in Sections \ref{sec:3D4}, \ref{sec:2E6}.

In general, the affine Hecke algebra $\ha(\mathbf{G}(k), \PH_k, \sigma)$ attached to some $\pi \in \mathrm{Irr}_{upt}(\mathbf{G}(k), \PH_k, \sigma)$ can have positive rank. A spectral transfer morphism from $\ha(\mathbf{G}(k), \PH_k, \sigma)$ to $\ha^{\IM}$ is given by a finite morphism $\mathcal{S}: \mathtt{Spec}\,Z(\ha) \to \mathtt{Spec}\, Z(\ha^{\IM})$ on the associated algebraic tori satisfying some extra conditions (see Def.~\ref{def:stm}). 

Among these technical conditions for spectral transfer morphism, the most important one is that a spectral transfer morphism should preserve the Plancherel measure in the following sense. We define a trace functional $\tau'$ on an affine Hecke algebra $\ha$. The spectral measure of this trace functional is called the \emph{Plancherel measure}. The density of this Plancherel measure is expressed by a function, called $\mu$-function, which is defined on the associated algebraic torus $T:=\mathtt{Spec} \, Z(\ha)$. A spectral transfer morphism between two pairs $(\ha_i, \tau'_i), i=1, 2$, is given by a finite morphism $\mathcal{S}: T_1 \to T_2$ which preserves the associated Plancherel measure. Roughly speaking, this means up to a non-zero rational constant, the pull-back of $\mu_2$ under $\mathcal{S}$ is equal to $\mu_1$. (See the condition (T3) in Def.~\ref{def:stm}.)

This paper aims at verifying that three maps are spectral transfer morphisms. In order to do this, we first recall the structure of the affine Hecke algebra $\ha(\mathbf{G}(k), \PH_k, \sigma)$ and define the $\mu$-function in Sections \ref{sec:genericAHA} and \ref{sec:mu} respectively. We will also recall residual cosets, which is necessary to give the definition of spectral transfer morphisms in Section \ref{sec:stm}. Then we turn to unipotent representations in Section \ref{sec:uptrep}. The materials here are borrowed from \cite{Opdl}, except for Sections \ref{sec:3D4} and \ref{sec:2E6}. We will recall the definition of unipotent representation and the associated affine Hecke algebra. In particular, we will give the parameters of the relevant affine Hecke algebras for classical groups in Section \ref{sec:para}. After these preparation, we give expressions of three morphisms in Section \ref{sec:stmexpressions}, and will prove in Section \ref{sec:pf} that these morphisms are spectral transfer morphisms.

\section{Affine Hecke Algebras and $\mu$-Functions}

\subsection{Generic affine Hecke algebras}\label{sec:genericAHA}

Let $\calr = (X, R_0, F_0, Y, R_0^\vee, F_0^\vee)$ be a reduced irreducible and semi-simple root datum, with a perfect pairing $\langle~,~\rangle : X \times Y \to \ZZ$. Here $F_0$ is a base of $R_0$. Let $\ZZ R_0$ be the root lattice span by the root system $R_0$. Then $\Omega_X:=X/\ZZ R_0$ is a finite abelian group. Let $\weylh_0=\weylh_0(R_0)$ be the finite Weyl group of the root system $R_0$, with a set $\mathtt{S}_0$ of distinguished generators, in bijection with $F_0$. We equip the vector space $E^\ast := \RR \otimes_{\ZZ} X$ with a Euclidean structure which is invariant under the action of $\weylh_0$. Following the convention in \cite{Opds}, we call elements of the set $R:=R_0^\vee \times \ZZ$ affine roots. If $a = \alpha^\vee + n \in R$, the reflection $r_a$ induced by $a$ is defined by $r_a(x)=x-a(x)\alpha$, where $a(x) = \langle x, \alpha^\vee \rangle + n$ for all $x\in E^\ast$. Such reflections are viewed as affine linear transformations on $E^\ast$.

The semi-direct product $\weylh:=X\rtimes \weylh_0$ is called the extended affine Weyl group. Inside $\weylh$ sits a normal subgroup $\weylh_a\simeq \ZZ R_0 \rtimes \weylh_0$ (the unextended affine Weyl group). $\weylh_a$ is an affine Coxeter group. Let $\tts_a$ be a set of distinguished generators of $\weylh_a$. 

Let $C \subset E^*$ be the fundamental alcove for the action of $\weylh_a$ on $E^\ast$. The closure $\overline{C}$ of $C$ is a fundamental domain for the action of $\weylh_a$ on $E^\ast$. The stabiliser 
$\mathrm{Stab}_{\weylh}(C) := \{ w \in \weylh: w(C)=C\}$ of $C$ in $\weylh$ is isomorphic to the quotient $\weylh/\weylh_a$, and hence isomorphic to $\Omega_X$.

On the affine Coxeter system $(\weylh_a, \tts_a)$ we have a canonical length function $l: \weylh_a \to \ZZ_+$ such that every distinguished generator $s \in \tts_a$ has length one. We extend $l$ to $\weylh$ by defining that $l(\omega)=0, \forall \omega \in \Omega_X$. Let $t_x \in \weylh$ denote the translation corresponding to $x\in X$. Then $l$ satisfies 
\begin{equation}
l(wt_x) = l(w) + \sum_{\alpha \in R_0^+ } \langle x, \alpha^\vee \rangle\end{equation}
for all $ w\in \weylh_0$ and all  $x \in X^+:= \{ x\in X: \langle x, \alpha^\vee \rangle \geq 0, \forall \alpha \in F_0 \}$.

Let $\CC[\mathbf{v}^\pm]$ be the Laurent polynomial ring in variables $\mathbf{v, v}^{-1}$. We define an associative algebra $\ha(\calr, m)$ over $\CC[\mathbf{v}^\pm]$, with distinguished $\CC[\mathbf{v}^\pm]$-basis $\{N_w: w \in \weylh \}$ parametrized by $w \in \weylh$, satisfying the following relations:
\begin{itemize}
\item [(i)] If $l(ww')=l(w)+l(w')$, then $N_w N_{w'} = N_{ww'}$.
\item [(ii)] $(N_s - \mathbf{v}^{m(s)}) (N_s + \mathbf{v}^{-m(s)})=0$ for all $s\in \tts_a$. 
\end{itemize}
Here the parameters $m(s) \in \ZZ$ are given by a function $m: \weylh\backslash \tts_a \to \ZZ$, which is defined on the $\weylh$-conjugacy classes of $\tts_a$. The algebra $\ha(\calr, m)$ is called the generic affine Hecke algebra (with parameter $m$). If there is no danger of ambiguity, we shall simply write $\ha(\calr)$ for $\ha(\calr, m)$.\\

\begin{rk}
\begin{itemize}
\item[(i)] Using the \emph{Bernstein presentation} (see e.g.~\cite[Thm.~2.6]{Opds}) for $\ha(\calr, m)$, one can show that the centre of $\ha(\calr, m)$ is $(\CC[\mathbf{v}^\pm][X])^{\weylh_0}$. 
\item[(ii)] From the definition, we note that $\ha(\calr, m)$ is not sensitive to the orientation of the double or triple edges in the Dynkin diagram of $R_0$ in case that the roots in $R_0$ have different lengths.
\end{itemize}
\end{rk}

Define $\tau: \ha(\calr) \to \CC[\mathbf{v}^\pm]$ by $\tau(\sum_w a_w N_w) = a_1$. Let $\ha_{v}(\calr)$ be the specialisation of $\ha(\calr)$ at $\mathbf{v}=v \in \RR_{>1}$. Then $\tau$ determines a family of $\CC$-valued linear functional $ \{ \tau_{v}\mid v \in \RR, v>1\}$ on the family of the algebras $\{\ha_v(\calr) \mid v \in \RR, v>1\}$. Define an anti-involution $\ast$ on $\ha_{v}(\calr)$ by $(\sum_w a_w N_w)^* = \sum_w \overline{a_w} N_{w^{-1}}$, then the functional $\tau_v$ is positive with respect to this anti-involution. Moreover, for all $v \in \RR_{>1}$, the triple $(\ha_{v}(\calr), \tau_v, *)$ form a type I Hilbert algebra in the sense of Dixmier \cite{Dixmier}.

The spectral measure associated with $\tau_v$ on $\ha_v$ is called the \emph{Plancherel measure}. Its density is given by the so-called $\mu$-function which we will review in next section.

\subsection{The $\mu$-function}\label{sec:mu}

Recall the root datum $\calr=(X, R_0, F_0, Y, R_0^\vee, F_0^\vee)$ underlying $\ha(\calr)$. Let $\mathcal{T}$ be the diagonalisable group scheme with character lattice $X \times \ZZ$. We view $\mathcal{T}$ as a group scheme over $\mathtt{Spec}\CC[\mathbf{v}^{\pm}] = \CC^\times$ via the homomorphism $\CC[\mathbf{v}^{\pm}] \to \CC[X \times \ZZ]$ defined by $\mathbf{v}^{n} \mapsto (0, n)$. The fibre at $v \in \CC^\times$ will be denoted by $\mathcal{T}_v$. The function $\mu(R_0, \mathbf{v})$ associated with $\ha(\calr, m)$ is a rational function on $\mathcal{T}$ with rational coefficient. 

To define $\mu(R_0, \mathbf{v})$ we need to introduce some notations. We define a function $\zm_R: R \to \ZZ$ by the rules that (i) $\zm_R$ is $\weylh$-invariant and (ii) $\zm_R (a) = m(s)$ if the simple reflection $s$ is conjugate to the reflection $r_{a}$ induced by $a \in R=R_0^\vee \times \ZZ$. This function $\zm_R$ gives rise to the functions $\zm_{\pm} : \weylh_0 \backslash R_0:\to \frac{1}{2}\ZZ$ defined as
\begin{equation}\label{eq:mplusminus}
\zm_+(\alpha) = \frac{\zm_R(\alpha^\vee) + \zm_R(1+\alpha^\vee)}{2}, \quad \zm_-(\alpha) = \frac{\zm_R(\alpha^\vee) - \zm_R(1+\alpha^\vee)}{2}
\end{equation}

One observes that $\zm_-(\cdot) = 0$ if and only if the two reflections $r_{\alpha^\vee}$ and $r_{1+\alpha^\vee}$ are $\weylh$-conjugate. So if $\zm_-(\cdot) \neq 0$, then $R_0$ must contain a component of type $\bbb$.

For each $\alpha \in R_0$, define
\begin{equation}
c_\zm(\alpha, \mathbf{v}): = \frac{(1 - \mathbf{v}^{-2\zm_+(\alpha)} \alpha^{-1})(1 + \mathbf{v}^{-2\zm_-(\alpha)} \alpha^{-1})}{1-\alpha^{-2}}.
\end{equation}
We put $c_\zm: = \prod_{\alpha \in R_0^+} c_\zm(\alpha, \mathbf{v})$. The Weyl group $\weylh_0$ acts on $c_\zm(\alpha)$, and hence on $c_\zm$, via its action on the roots. The $\mu$-function of $\ha(\calr, m)$ is defined as 
\begin{equation}
\mu(R_0, \mathbf{v}) = \mathbf{v}^{-2\zm_{\weylh}(w_0)} \frac{d(\mathbf{v})}{(w_0. c_\zm) c_\zm} = \frac{d(\mathbf{v})}{\mathbf{v}^{2\zm_{\weylh}(w_0)}} \prod_{\alpha \in R_0} \frac{1}{c_\zm(\alpha, \mathbf{v})}
\end{equation}
where $w_0 \in \weylh_0$ is the longest element, $\zm_{\weylh}: \weylh \to \ZZ$ is defined by 
\begin{equation*}
\zm_{\weylh}(w) = \sum_{a \in R^+ \cap w^{-1}(R^-)} \zm_R (a),
\end{equation*}
and the factor $d(\mathbf{v})$ is a rational function in variable $\mathbf{v}$. Clearly, the $\mu$-function is $\weylh_0$-invariant.

The explicit expression of $d(\mathbf{v})$ will be given in Section 3. Note that $\mu(\emptyset, \mathbf{v})=d(\mathbf{v})$. Here we point out that we shall regard $d(\mathbf{v})$ as a normalisation factor in the following sense: let $1 \in \weylh$ be the identity element of the extended affine Weyl group, then the functional $N_w \mapsto \delta_{w, 1} d(\mathbf{v})$ defines a trace $\tau'(\cdot):=d(\mathbf{v})\tau(\cdot)$ on $\ha(\calr, m)$. 

We call $\mu_{\zm, d}(R_0, v):= \mu(R_0, \mathbf{v})|_{\mathbf{v}=v}$ the $\mu$-function associated with $(\ha_v(\calr, m), d)$, where $d:=d(\mathbf{v})|_{\mathbf{v}=v}$.

\subsection{Residual cosets and residual points}\label{sec:respt}

Let $T=\mathcal{T}_v$ be the complex torus defined as in Section \ref{sec:mu}, so that the character lattice of $T$ is $X$. Let $L \subset T$ be a coset of a subtorus $T^L$ (say) of $T$. We define (where $\epsilon =\pm$ is a sign)
\begin{equation}\label{def:respt}
 p_{\epsilon}(L) = \{\alpha \in R_0: \epsilon \alpha|_{L} = v^{-2\zm_{\epsilon}(\alpha)} \}, \quad z_{\epsilon}(L) = \{\alpha \in R_0: \epsilon \alpha|_{L} = 1 \}
\end{equation}
The coset $L \subset T$ is said to be \emph{residual} if 
\begin{equation}
|p_+(L) | + |p_-(L) | - |z_+(L) |-|z_+(L) | = \text{codim}(L). 
\end{equation}

Suppose $L$ is a residual coset for $\mu_{m, d}(R_0, v)$ with $p_{\epsilon}(L)$ and $z_{\epsilon}(L)$ as defined above, by the \emph{regularisation of $\mu$ along $L$} we mean the following expression
\begin{equation}\label{eq:regularisation}
\mu_{\zm, d}^L(R_0, v): = \frac{(d/v^{2\zm_{\weylh}(w_0)}) \prod_{\alpha \in R_0 \backslash z_-(L)} (1+\alpha^{-1}) \prod_{\alpha \in R_0 \backslash z_+(L)} (1-\alpha^{-1})}
{\prod_{\alpha \in R_0 \backslash p_-(L)} (1+v^{-2\zm_-(\alpha)} \alpha^{-1}) \prod_{\alpha \in R_0 \backslash p_+(L)} (1-v^{-2\zm_+(\alpha)} \alpha^{-1}) }.
\end{equation} 

The regularisation $\mu_{\zm, d}^L(R_0, v)$ restricting onto $L$ defines a nonzero rational function $\mu_{\zm, d}^{(L)}(R_0, v)$ in $v$ over the residual coset $L$. In particular, if $L=\{ \respt \}$ is just a residual point, we call the expression $\mu_{\zm, d}^{(\{\respt \}) }(R_0, v)$ the ``residue'' of the $\mu$-function at the residual point $\respt$, and we shall use the notation $\Res(\mu_{\zm, d}(R_0, v), \respt)$ for $\mu_{\zm, d}^{(\{\respt \}) }(R_0, v)$, to distinguish the role of the residual point $\respt$. (Compare with the computation of residue of a one-variable meromorphic function at a simple pole.) 

Since $T$ is a complex torus, it has a polar decomposition $T=T_u T_+$, where $T_u$ is unitary and $T_+$ is real split. For $\respt \in T$, we usually write $s$ for the unitary part of $\respt$.

\section{The Spectral Transfer Morphisms for Affine Hecke Algebras} \label{sec:stm}


Given two normalised affine Hecke algebras $(\ha(\calr_1, m_1), d_1)$ and $(\ha(\calr_2, m_2), d_2)$, let $T_1, T_2$ be the associated tori. Let $L=\respt T^L \subset T_2$ be a residual coset, with $\respt \in T_L$ a residual point (here $T^L$ is a subtorus of $T$). Let $\weylh_{2, 0}$  be the Weyl group associated with the root datum of $\calr_2$. Denote by $\weylh_{2, 0}(L)$ the quotient group of the stabiliser of $L$ in $\weylh_{2, 0}$ by the pointwise stabiliser. \\

\begin{definition}\cite[Def.~5.1]{Opds}\label{def:stm}
Consider a morphism $\phi_T: T_1 \to T_2$, whose range is a $L$. We call $\phi_T$ a \emph{spectral transfer morphism} if it satisfies the following conditions:
\begin{enumerate}
\item[(T1)] $\phi_T$ is a finite morphism (in the sense of algebraic geometry);
\item[(T2)] $\phi_T$ maps the identity $e$ of $T_1$ into $L\cap T_L$, and if we declare $\phi_T(e)$ to be the base point of $L$, then $\phi_T$ is a homomorphism of algebraic tori.
\item[(T3)] There exists a constant $c \in \QQ^\times$ such that $\phi_T^* \bigg(\mu^{(L)}_{\zm_2, d_2}(\calr_2) \bigg) = c \mu_{\zm_1, d_1}(\calr_1)$.\\
\end{enumerate}
\end{definition}

\begin{rk}
There is an extra technical condition (T4). However, to deal with unipotent representations we do not need it and hence we do not state it. In practice, the conditions (T1), (T2) and (T4) are usually easily verified. So the key point is the verification of (T3). (See also \cite[Def.~5.1]{Opds} and \cite[Def.~4.3.1, Remark 4.3.2]{OpdBernstein}.)
\end{rk}

A surjective spectral transfer morphism is called a \emph{spectral covering}. Examples of spectral coverings can be found in \cite[7.1.3]{Opds}. However, the main examples in this paper all satisfy $\dim T_1 < \dim T_2$, and hence are not surjective. \\

\begin{exa} \label{exa:Weylgpelement}
Let $\weylh_0$ be the finite Weyl group associated to the root datum $\calr$. Every $w \in \weylh_0$ acts on the associate complex torus as automorphism. Since the $\mu$-function is $\weylh_0$-invariant, this morphism clearly satisfies the condition (T3).

There are other kinds of spectral transfer morphisms which are induced from the algebra structure of $\ha_i$. See \cite[7.1.1]{Opds}.  
\end{exa}

\section{Unipotent Representations and Affine Hecke Algebras}\label{sec:uptrep}

\subsection{Unipotent representations}

We will define unipotent representations for $\mathbf{G}(k)$ and discuss the affine Hecke algebras derived from unipotent representations. Firstly, let us give the assumptions on the ground field $k$ and on the linear algebraic group $\mathbf{G}$. 

Let $k$ be a non-archimedean local field with finite residue field $\mathbf{F}_q$. Fix a separable closure $k_s$ of $k$. Let $\KK \subset k_s$ be the maximal unramified extension of $k$, with residual field $\overline{\FF}_q$, an algebraic closure of $\FF_q$. Recall that there exist isomorphisms of Galois groups: $\mathrm{Gal}(\KK/k) \simeq \mathrm{Gal}(\overline{\FF}_q/\FF_q) \simeq \hat{\ZZ}$. The \emph{geometric Frobenius element} $\mathrm{Frob}$, whose \emph{inverse} induces the automorphism $x \mapsto x^{q}$ for any $x \in \overline{\mathfrak{F}}$, is a topological generator of $\mathrm{Gal}(\KK/k)$.

We assume that $\mathbf{G}$ is a connected absolutely simple algebraic group of adjoint type defined and quasi-split over $k$ such that $\mathbf{G}$ splits over $\KK$. The Galois group $\mathrm{Gal}(k_s/k)$ acts on $\mathbf{G}(k_s)$. Since $\mathbf{G}$ is $\KK$-split, this action factors through an action of $\mathrm{Gal}(\KK/k)$ on $G:=\mathbf{G}(\KK)$. It turns out that the action of $\mathrm{Gal}(\KK/k)$ on $G$ is completely determined by the action of $\Frob$ on it. We shall denote $F \in \mathrm{Aut}(G)$ the automorphism induced by $\Frob$.

Recall that the equivalence classes of inner forms of $\mathbf{G}$ are parameterized by the Galois cohomology set $H^1(k, \mathbf{G}_{ad}) = H^1(k, \mathbf{G})$. By a theorem of Steinberg (cf.~\cite[Thm.~1.9]{St}), we know that $H^1(\mathrm{Gal}(k_s/\KK), \mathbf{G})$ is trivial. We thus have an isomorphism $H^1(k, \mathbf{G}) \cong H^1(F, G):=H^1(\mathrm{Gal}(\KK/k), G)$. 
Let $z\in Z^1(F, G)$ be a cocycle and denote $u:=z(\Frob) \in G$. Let $F_u:=\mathrm{Ad}(u) \circ F$ be an inner twist of $F$ by $u$. Then $F_u$-action on $G$ defines a $k$-structure on $G$. If $\mathbf{G}^u$ is the inner form of $\mathbf{G}$ corresponding to $u$, then $\mathbf{G}^u(k)= G^{F_u}$. In particular, when $[z]=1 \in H^1(F, G)$ is the privileged element, we have $\mathbf{G}(k) = G^F$. 

We can find a parahoric subgroup $\PH \subset G$ such that $F_u(\PH)=\PH$. Let $\PH_+$ be the pro-unipotent radical of $\PH$. The quotient $\PH/\PH_+$ is reductive. Since $\PH_+$ is connected, we can apply the pro-algebraic version of Lang's theorem and obtain $(\PH/\PH_+)^{F_u} = \PH^{F_u}/\PH_+^{F_u}$, which is a finite reductive group. We shall denote this finite group by $\overline{\PH}^{F_u}$. 

A smooth representation $\pi$ of $G^{F_u}$ is called \emph{unipotent}, if there exists an $F_u$-stable parahoric subgroup $\PH$ of $G$, with a cuspidal unipotent representation $\sigma$ of $\overline{\PH}^{F_u}$, such that the $\PH_+^{F_u}$-invariants of $\pi$ contain $\sigma$. (Here we lift $\sigma$ to a representation of $\PH^{F_u}$ via the natural projection $\PH^{F_u} \to \overline{\PH}^{F_u}$). The pair $\utype:=(\PH^{F_u}, \sigma)$ a \textit{type} of $G^{F_u}$ in the sense of Bushnell--Kutzko \cite{BK}. \\

\begin{definition}
Let $\PH \subset G$ be a $F_u$-stable parabolic subgroup, and let $\sigma$ be a cuspidal unipotent representation of $\overline{\PH}^{F_u}$. We call the pair $\utype=(\PH^{F_u}, \sigma)$ a \emph{unipotent type} of $G^{F_u}$, and denote by $\mathrm{Irr}(G^{F_u}, \utype)$ the totality of isomorphism classes of irreducible unipotent representations of $G^{F_u}$ attached to a unipotent type $\utype=(\PH^{F_u}, \sigma)$ of $G^{F_u}$. 
\end{definition}

\subsection{Unipotent affine Hecke algebras}\label{sec:AHAupt}


We fix a unipotent type $\utype=(\PH^{F_u}, \sigma$) of $G^{F_u}$ in this section. Define a representation $\cind_{\PH^{F_u}}^{G^{F_u}} \sigma$ of $G^{F_u}$. Let $V$ be the vector space of this representation. The endomorphism algebra $\End(V)$ is isomorphic to the $(\PH^{F_u}, \sigma)$-spherical Hecke algebra of $G^{F_u}$. We denote this $(\PH^{F}, \sigma)$-spherical Hecke algebra by $\ha_{\utype}$. Let $\mathrm{Irr} \left( \ha_{\utype}\right)$ be the set of isomorphism classes of (finite dimensional) simple $\ha_{\utype}$-modules. From the work of Moy--Prasad \cite{MP}, Morris \cite{Morris} and Lusztig\cite{Lusztig}, we know that\\
\begin{prop}\label{prop:repandhmod}
\begin{itemize}
\item[(i)] There is a natural bijection between the sets $\mathrm{Irr}\left(G^{F_u}, \utype \right)$ and $\mathrm{Irr}\left(\ha_{\utype}\right)$.
\item[(ii)] Two sets $\mathrm{Irr}(G^F, \utype)$ and $\mathrm{Irr}(G^F, \utype')$ are either disjoint or identical. They are equal if and only if there exists an element $g\in G^F$ which conjugates $\PH$ to $\PH'$, and $\sigma$ to a representation isomorphic to $\sigma'$.
\end{itemize}
\end{prop}

The structure of $\ha_{\utype}$ has been determined by Lusztig \cite{Lusztig}. Briefly, $\ha_{\utype}$ is an extension of an affine Hecke algebra by a finite group. This finite group is commutative, and can be described in terms of the group of diagram automorphisms of the affine Dynkin diagram $\dyn(\Sigma_0^{(1)})$ of $G$. Here $\Sigma_0^{(1)}$ is the affine root system attached to $\mathbf{G}$. 


Let $\mathcal{I}\subset \PH$ be an Iwahori subgroup of $G$. Define a group $\Omega:=N_G(\mathcal{I})/\mathcal{I}$. It is known that $\Omega$ is finite and abelian (cf.~\cite{DeRe}). After conjugation (if necessary), we may and will assume that $\mathcal{I}$ is standard. Thus, the group $\Omega$ acts on $\dyn(\Sigma_0^{(1)})$ of $\mathbf{G}$ as diagram automorphisms. 

Let $\tisfI$ denote the set of all nodes of $\dyn(\Sigma_0^{(1)})$. Suppose that $\PH$ in the type $\utype=(\PH^{F_u}, \sigma)$ corresponds to a subset $\mathsf{J} \subset \tisfI$. The Galois group $\mathrm{Gal}(k_s/k)$ acts on $\dyn(\Sigma_0^{(1)})$. This action factors through the action of $\mathrm{Gal}(\KK/k)$, because $\mathbf{G}$ is $\KK$-split. Denote by $\theta$ the diagram automorphism of $\dyn(\Sigma_0^{(1)})$ induced by $F$ (the automorphism of $G$ induced by $\Frob$). We also have the twisted group $\theta_u$, induced by $F_u$. Since $\PH$ is $F_u$-stable, we see that $\theta_u$ preserves $\mathsf{J}$.

Furthermore, let $\Omega^{\theta_u}$ denotes the elements in $\Omega$ commuting with $\theta_u$. Since $\Omega$ is abelian, we thus have $\Omega^{\theta_u}=\Omega^\theta$. The action of $\Omega$ on $\dyn(\Sigma_0^{(1)})$ induces an action of $\Omega^\theta$ on the $\theta_u$-orbits of this Dynkin diagram. Denote by $\Omega^{\theta}(\PH)$ the isotropy group of $\PH$. Then $\Omega^\theta(\PH)$ acts on the set $(\tisfI - \mathsf{J})/\theta_u$ of $\theta_u$-orbits as permutations. Let $\Omega_1^{\theta}(\PH)$ be the pointwise stabilizer of $(\tisfI-\mathsf{J})/\theta_u$, and let $\Omega_2^{\theta}(\PH):= \Omega^{\theta}(\PH)/\Omega_1^{\theta}(\PH) $. 

Lusztig \cite[1.18, 1.19, 1.20]{Lusztig} described the structure of $\ha_{\utype}$. It decomposes into a direct sum decomposition of two-sided ideals:
$$
\ha_{\utype}= \ha^\psi \rtimes \Omega_1^{\theta}(\PH)
$$ 
such that: every ideal $\ha^\psi$ is isomorphic to the other, and they are parameterised by irreducible characters $\psi$ of $\Omega_1^{\theta}(\PH)$. He then constructed an isomorphism between $\ha^\psi$ and an affine Hecke algebra obtained as a generic affine Hecke algebra specialised at $\mathbf{v}=q^{1/2}$, and thus gave the Iwahori--Matsumoto presentation of $\ha^\psi$ for all types of simple algebraic groups which are inner forms of split groups. Associated to $\ha^\psi$ is an affine Weyl group $W_{\utype}$ with a set $S_{\utype}$ of distinguished generators. The structure of $W_{\utype}$ as well as the parameter function $m_{\utype}$ can be entirely determined by the data of $\mathbf{G}$ and $\utype$ (cf.~\cite{Morristype}). However, there is no canonical way to choose the set $S_{\utype}$. But different choices of $S_{\utype}$ give us algebras which are essentially the same, in the sense that there exist admissible isomorphisms (cf.~\cite[2.1.7]{Opds} and \cite[text around Eq.~(19)]{Opdl}) among them. 

The following proposition is a consequence of the theory of \textit{types} due to \cite{BK}. It is an enhancement of the bijection (i) in Proposition \ref{prop:repandhmod}. The claim of bijection is proved at the end of \cite{MP} as an equivalence of two categories of representations. Using the main theorem in \cite{Opds} we obtain the other claim.\\ 

\begin{prop}
Let $\utype=(\PH^{F_u}, \sigma)$ be a unipotent type of $G^{F_u}$. The bijection between $\mathrm{Irr} (G^{F_u}, \utype)$ and $\mathrm{Irr} (\ha_{\utype})$ respects the notion of temperedness, and preserves the Plancherel measure on the level of irreducible tempered representations. 
\end{prop}

We call $\ha^\psi$ a unipotent affine Hecke algebra, and denote it by $\ha_{upt}$. In particular, we can define discrete series representations for unipotent affine Hecke algebras, and study their formal degrees in place of the formal degrees of unipotent discrete series representations of $G^F$.


In Section \ref{sec:para}, we will discuss the parameters for $\ha_{upt}$ when $\mathbf{G}$ is isogenous to a classical group (except for the special linear group).
   
\subsection{Cuspidal unipotent characters and the normalisation factor $d$}


From \cite[\S 13.7]{Carter} we know that cuspidal unipotent characters are rare for classical finite groups of Lie type. The root system type must verify selective conditions. For instance, an odd orthogonal group has one cuspidal unipotent character only if the rank is of the form $s(s+1)$ for some positive integer $s$, and has no such character otherwise. Therefore, we can determine the subset $\mathsf{J} \subset \tisfI$ corresponding to $\PH$ such that $\overline{\PH}^{F_u}$ has a cuspidal unipotent character. Consequently, we are able to explicitly define the normalisation factor $d(\mathbf{v})$ of the $\mu$-function of $\ha_{upt}$. We proceed as follows.

We normalise the Haar measure on $\mathbf{G}(k)$ as \cite{DeRe}, then the volume of a parahoric subgroup $\PH^{F_u}$ equals to 
\begin{equation}\label{eq:Haarmeasure}
\mathrm{Vol}(\PH^{F_u})= v^{-a} |\overline{\PH^{F_u}}|,
\end{equation}
where $v$ is the positive square root of $q$, $a$ is the rank of $\mathbf{G}$ regarded as an algebraic group defined over a separable closure $k_s$ of $k$, and $\overline{\PH_k}$ is the the reductive quotient of $\PH_k$ by its pro-unipotent radical $\PH_{k, +}$. So $\overline{\PH_k}$ is a finite group of Lie type. Using the results in \cite[\S 2.9]{Carter} we can determine $|\overline{\PH_k}|$, the order of $\overline{\PH_k}$. 

For a unipotent discrete series representation belongs to $\mathrm{Irr}_{upt}(G^{F_u}, \utype)$ with $\utype=(\PH^{F_u}, \sigma)$,  via the analysis in \cite[\S 2.4.1]{Opdl}, the normalisation factor $d=d(v)$ is given by \cite[Eq.~(24)]{Opdl}
\begin{equation}
d(v)=|\Omega_1^\theta(\PH)|^{-1} \mathrm{Vol}(\PH^{F_u})^{-1} \deg (\sigma).
\end{equation}
One can check case-by-case that $d(v)$ satisfies that $d(v)=d(v^{-1})$. 

\subsection{The case of non-split but quasi-split exceptional groups}

Let $\mathbf{G}$ be an adjoint $k$-group of type ${}^3\ddd_4$ or ${}^2\eee_6$. In this section, we show that for a given supercuspidal unipotent representation $\pi$ of $\mathbf{G}(k)$, there is only Weyl group orbit $\weylh_0 \respt$ of residual point verifying the following equation 
\begin{equation}\label{eq:fdeg=resmu}
\fdeg(\pi, q) = C_\pi \Res(\mu^{\IM}, \respt)
\end{equation}
for some constant $C_\pi \in \QQ^\times$ independent of $q$. In other words, the spectral transfer morphism from $\ha_{upt}$ (in this case $\ha_{upt}=\CC$) to the Iwahori--Hecke algebra $\ha^{\IM}$ is well-defined.

\subsubsection{The group of type ${}^3\ddd_4$}\label{sec:3D4}

Let $\mathbf{G}$ be an adjoint $k$-group of type ${}^3\ddd_4$. Under the Frobenius action, the relative root system is $\mathrm{G}_2$. So the underlying root system of the unipotent affine Hecke algebra is $R_0=\mathrm{G}_2$. Let us denote by $F_0=\{\alpha_1, \alpha_2\}$ the base of $R_0$, where $\alpha_1$ is a short root and $\alpha_2$ is a long root. In terms of \cite[Eq.~(4) and (8)]{Opds}, the parameters of the associated $\mu$-function are given as follows: $m_+(\alpha_1)=1, m_+(\alpha_2)=3$. The associated Iwahoir-Hecke algebra is denoted as $\mathrm{G}_2(3, 1)[q]$ in notations of \cite[3.2.1]{Opdl}, and the normalisation faction $d^{\IM}$ in the $\mu^{\IM}$-function is given explicitly as 
$$
d^{\IM}= \frac{v-v^{-1}}{v^3-v^{-3}} \frac{1}{(v-v^{-1})^{2}} = \frac{1}{(v-v^{-1}) (v^3-v^{-3})} = \frac{q^2}{(q-1)(q^3-1)}.
$$ 

From \cite[\S 13.7]{Carter} we know that there are in total two cuspidal unipotent characters of the finite reductive group of type ${}^3\ddd_4$. They are indicated in the right-most column at Table \ref{tab:3D4}. The algebra $\mathrm{G}_2(3, 1)[q]$ has four orbits of residual orbits. Among them, we find two residual orbits with representatives $\respt_1 = (1, q)$ and $\respt_2=(-q^{-3}, q)$ verify Eq.~\eqref{eq:fdeg=resmu} for these two cuspidal unipotent characters respectively. These are the only residual orbits of $\mathrm{G}_2(3, 1)[q]$ supporting supercuspidal unipotent representations. 


\subsubsection{The group of type ${}^2\eee_6$}\label{sec:2E6}

Now let $\mathbf{G}$ be an adjoint $k$-group of type ${}^2\eee_6$. Under the Frobenius action we have a relative system $\mathrm{F}_4$. So the root system of the unipotent affine Hecke algebra is $R_0= \mathrm{F}_4$. The parameter function $m_+$ takes value $2$ (resp.~$1$) at long (resp.~short) roots in $R_0 = \mathrm{F}_4$. The Iwahori--Hecke algebra is denoted as $\mathrm{F}_4(2, 1)[q]$.

There are 24 positive roots in $\mathrm{F}_4$. The normalisation factor of the $\mu^{\IM}$-function is
$$
d= (v-v^{-1})^{-2} (v^2-v^{-2})^{-2} = q^{-3} (q-1)^{-2} (q^2-1)^{-2},
$$
and there are eight residual orbits. From \cite[\S 13.7]{Carter} we know that there are in total three cuspidal unipotent characters of the finite reductive group of type ${}^2\eee_6$. They are indicated in the right-most column at Table \ref{table:2E6}. Note that two cuspidal unipotent characters ${}^2\eee_6[\theta]$ and ${}^2\eee_6[\theta^2]$ have the same degree. 

With the help of Maple, we obtain Table \ref{table:2E6} below. We see that two orbits with representatives $\respt_1=(1,1,q,1)$ and $\respt_2=(q^2, \theta q^{-4}, q, q)$ (where $\theta$ is a primitive third root of unity) verify Eq.~\eqref{eq:fdeg=resmu} for these three cuspidal unipotent characters respectively. These are the only residual orbits of $\mathrm{F}_4(2, 1)[q]$ supporting supercuspidal unipotent representations. 

\begin{rk}
The term $s$ at the left-most column of these tables has two-fold meaning. Firstly, it is exactly the semi-simple part of the image of the unramified local Langlands parameter associated with the supercuspidal unipotent representation $\pi$ (cf.~Eq.~\eqref{eq:imagephi}). It is known that $s$ must be a torsion element in $G^\vee$, and hence $s$ corresponds to a node in the Kac diagram (cf.~\cite[\S 3.8 and Table 1]{Retorsion}). We indicate $s$ by the type of the complement of this node.) Secondly, it is the compact part of the residual point $\respt \in T$ under the polar decomposition of the complex torus $T$. The bijection mentioned in \cite[Prop.~B.1]{Opd3} shows that these two interpretations of $s$ are compatible. 
\end{rk}

\begin{table}[htbp]
  \caption{Cuspidal unipotent characters of ${}^3\ddd_4$.}
   \label{tab:3D4}
   \setlength\extrarowheight{2pt}
   \centering
    \begin{tabular}{|p{3em}|p{4em}|p{10em}|p{4em}|}
   \hline
   $s$ & $\respt$ & $\Res(\mu_{q'}^{\IM}, \respt)$ & $\mathcal{U}^0$ \\
      \hline
  $\mathrm{G}_2$ & $(q^3, q)$ & $(\Phi_5 \Phi_9)/(\Phi_2^2 \Phi_3^2 \Phi_6^2 \Phi_{12})$&\\ [1.5ex]
    \hline  
 $\mathrm{G}_2(a_1)$ & $(1,q)$ & $(\Phi_2^2 \Phi_3^2 \Phi_6^{2})^{-1}$ & ${}^3\ddd_4[1]$\\[1.2ex]
 \hline
 $\aaa_2$ & $(q^3, \vartheta q^{-3})$ & $\Phi_3^{-2} \Phi_6^{-2} \Phi_9 \Phi_{12}^{-1}$ & \\[1.2ex]
 \hline
 $\aaa_1 \aaa_1$ & $(-q^{-3}, q)$ & $(\Phi_2^2 \Phi_6^2 \Phi_{12})^{-1}$ & ${}^3\ddd_4[-1]$ \\[1.5ex]
 \hline    
\end{tabular}
\end{table}

\begin{table}[htbp]
  \caption{Cuspidal unipotent characters of ${}^2\mathrm{E}_6$.}
   \label{table:2E6}
   \setlength\extrarowheight{2.8pt}
   \centering
    \begin{tabular}{|p{3em}|p{9em}|p{13em}|p{7em}|}
   \hline
   $s$ & $\respt$ & $\Res(\mu_{q'}^{\IM}, \respt)$ & $\mathcal{U}^0$ \\
      \hline
  $\mathrm{F}_4$ & $(q^2, q^2, q, q)$ & $\frac{\Phi_5 \Phi_7 \Phi_{11} \Phi_{16}}{\Phi_2^6 \Phi_3^2 \Phi_4^2 \Phi_6^3 \Phi_{10} \Phi_{12} \Phi_{18}}$ &\\ [1.7ex]
  $\mathrm{F}_4(a_1) $ & $(q^2, 1, q, q)$ & $\Phi_2^{-6} \Phi_3^{-2} \Phi_4^{-2} \Phi_5^{2} \Phi_6^{-3} \Phi_7 \Phi_{12}^{-1}$  & \\ [1.4ex] 
  $\mathrm{F}_4(a_2)$ & $(q^2, 1, q, 1)$ & $\Phi_2^{-6} \Phi_3^{-2} \Phi_4^{2} \Phi_5 \Phi_6^{-3} \Phi_{10}^{-1}$ & \\ [1.4ex]
  \hline
  $\mathrm{F}_4(a_3)$ & $(1, 1, q, 1)$ & $(\Phi_2^{6} \Phi_3^{2} \Phi_4^{2} \Phi_6^3)^{-1}$ & ${}^2\mathrm{E}_6[1]$ \\ [1.3ex]    
 \hline    
 $\aaa_1 \ccc_3$ & $(-q^{-7}, q^2, q, q)$ & $\Phi_{2}^{-6} \Phi_5 \Phi_{6}^{-3} \Phi_{18}^{-1}$ & \\ [1.1ex]
 \hline
 $\aaa_2 \aaa_2$ & $(q^2, \theta q^{-4}, q, q)$ & $(\Phi_{3}^{2}\Phi_{6}^{3}\Phi_{12}\Phi_{18})^{-1}$ & ${}^2\mathrm{E}_6[\theta]$, ${}^2\mathrm{E}_6[\theta^2]$ \\ [1.3ex]
 \hline
 $\aaa_3 \aaa_1$ & $(q^2, q^2, \sqrt{-1} q^{7/2}, q)$ & $\Phi_{2}^{-6}\Phi_{4} \Phi_{6}^{-3}\Phi_{8}\Phi_{10}^{-1}\Phi_{18}^{-1}$ & \\ [1.2ex]
 \hline
 $\bbb_4(1)$ & $(q^2, q^2, q, -q^{-8})$ & $\Phi_{2}^{-6} \Phi_4^{-2} \Phi_{5}\Phi_{6}^{-3}\Phi_{7}\Phi_{12}^{-1}\Phi_{16}\Phi_{18}^{-1}$ & \\ [1.2ex]
 \hline
 $\bbb_4(2)$ & $(q^2, 1, q, -q^{-5})$ & $\Phi_{2}^{-6} \Phi_{4}^{-2} \Phi_{5} \Phi_{6}^{-3} \Phi_{8}\Phi_{10}^{-1}\Phi_{12}^{-1}$ & \\[1.2ex]
 \hline
\end{tabular}
\end{table}

\subsection{The case of classical groups}\label{sec:para}

We still need to determine the parameters $\zm_\pm(\alpha)$ for each $\weylh_0$-conjugacy classes of roots $\alpha \in R_0$. For classical groups, these parameters can be given via another two closely related parameters $m_\pm$ such that $4m_\pm \in \ZZ$ and $2(\emm \pm \emp) \in \ZZ$. These parameters $m_\pm$ are determined from two natural numbers $a$ and $b$ coming from the condition that $\overline{\PH}$ has a cuspidal unipotent character. The relationship between $\{a, b\}$ and $\{\emm, \emp\}$ (here we understand them as multi-sets, i.e.~the elements can be equal) are explicitly given by \cite[3.2.1]{Opdl}. An important fact is, the data $\{a, b\}$ come from the representation (or arithmetic) side, while the data $\{m_-, m_+\}$ come from the spectral (or geometry) side, and the correspondence of them should be viewed as a correspondence satisfying Langlands' philosophy. 

Let $\mathcal{V}$ be the parameter space of the pair $\{\emm, \emp\}$. Recall that we have $4m_\pm \in \ZZ$. We partition $\mathcal{V}$ into six disjoint subsets as follows. 
  \begin{equation}
   \begin{array}{ll} 
    \{\emm, \emp\} \in \parav^{\rmone} & \text{ iff } m_{\pm} \in \ZZ/2 \text{ and } \emm - \emp \notin \ZZ, \\
     \{\emm, \emp\}  \in \parav^{\rmtwo} & \text{ iff } m_{\pm} \in \ZZ+1/2 \text{ and } \emm - \emp \in \ZZ, \\
      \{\emm, \emp\}  \in \parav^{\rmthree} & \text{ iff } m_{\pm} \in \ZZ \text{ and } \emm - \emp \notin 2\ZZ, \\
     \{\emm, \emp\}  \in \parav^{\rmfour} & \text{ iff } m_{\pm} \in \ZZ \text{ and } \emm - \emp \in 2\ZZ, \\
     \{\emm, \emp\}  \in \parav^{\rmfive} & \text{ iff } m_{\pm} \in \ZZ \pm 1/4 \text{ and } \dem \neq \dep, \\
      \{\emm, \emp\} \in \parav^{\rmsix} & \text{ iff } m_{\pm} \in \ZZ \pm 1/4 \text{ and } \dem = \dep.
   \end{array}
   \end{equation} 
The numbers $\delta_\pm$ appearing in cases V and VI are determined by the following method. If $m \in \ZZ \pm 1/4$, we write 
\begin{equation}\label{eq:kappa}
|m| = \kappa + \frac{2\epsilon - 1}{4},
\end{equation}
where $\epsilon \in \{0, 1\}$ and $\kappa \in \ZZ_{\geq 0}$. The integer $\kappa$ determines a number $\delta \in \{0, 1\}$ by the rule that $\kappa - \delta \in 2\ZZ$. In other words, $\delta$ indicates the parity of $\kappa$.\\

\begin{rk}
The isogeny type of the corresponding classical groups are: unitary groups in case I; odd special orthogonal groups in case II; symplectic groups in cases III and V; even special orthogonal groups in cases IV and VI. Cases V and VI correspond to the case that the centre of the spin group on the Langlands dual side acting non-trivially on the enhanced Langlands parameters. Following \cite{Opdl}, we call them the \emph{extra-special} cases.
\end{rk}

Finally, the parameters $\zm_\pm(\alpha)$ are obtained from $m_\pm$ by the following rules: Observe that from Lusztig's list \cite[\S 7]{Lusztig} we can assume that the relevant root system $R_0$ of the Iwahori--Hecke algebra always have type $\bbb$. Let $(t_i)_{1\leq i \leq l}$ be the standard basis of the character lattice which contains $\bbb_l$. We write the roots multiplicatively. Then:
\begin{equation}
\zm_-(t_i t_j ^{\pm} )=0, \;\zm_+(t_i t_j^\pm) = \frb;\; \zm_-(t_i) = \frb \emm, \;\zm_+(t_i) = \frb \emp,
\end{equation}
where $\frb=1$ if both $\emp \pm \emm \in \ZZ$, and $\frb=2$ else. 


\subsection{Three morphisms}\label{sec:stmexpressions}

We now look at three morphisms which are from unipotent affine Hecke algebras associated with \emph{classical} groups. Opdam defined them in \cite[3.2.6]{Opdl} and claimed (without proof) there that these morphisms are spectral transfer morphisms. We now give their definitions and prove that they are indeed spectral transfer morphisms. In what follows, we will write $T^n$ to indicate that the rank of the algebraic torus $T$ is $n$.

\begin{enumerate}
\item[(1)] For $m_{\pm} \in \ZZ + 1/2$, define a homomorphism $\phi_{T, \emm, \emp}: T^n \to T^{n+\emp-1/2}$ of algebraic tori over $\CC$ by 
\begin{equation}\label{eq:stmphi}
\phi_{T, \emm, \emp} (t_1, \ldots, t_n) = (t_1, \ldots, t_n, v^{\frb}, v^{3\frb}, \ldots, v^{2\frb(\emp-1)}).
\end{equation}

\item[(2)] For $\emp \in \ZZ_{>0}$, define a homomorphism $\psi_{T, \emm, \emp}: T^n \to T^{n+2\emp -2}$ of algebraic tori over $\CC$ by 
\begin{equation}\label{eq:stmpsi}
\begin{aligned}
&\psi_{T, \emm, \emp}(t_1, \ldots, t_n) \\
= &(t_1, \ldots, t_n, 1, q^{\frb}, q^{\frb}, q^{2\frb}, q^{2\frb}, \ldots, q^{\frb(\emp-2)}, q^{\frb(\emp-2)}, q^{\frb(\emp-1)}).
\end{aligned}
\end{equation}

\item[(3)] (The extra-special cases $\parav^{ \rmfive}$ and $\parav^{\rmsix}$). For $m_{\pm} >0$, set 
\begin{equation}\label{eq:l}
l:= \begin{cases}
2n+(a/2)(a+1)+2b(b+1) &\text{ if } \parax = \rmfive, \\
2n+(a/2)(a+1)+2b^2-\dep & \text{ if } \parax = \rmsix.
\end{cases}
\end{equation} 
Let $\kappa_{\pm}, \epsilon_{\pm}$ be defined as in (\ref{eq:kappa}). If we let $l_{\pm} := \kappa_{\pm}(\kappa_{\pm}+\epsilon_{\pm}-1/2)$, then $l=2n+\lfloor l_- \rfloor + \lfloor l_+ \rfloor$ in both parameter types $\rmfive, \rmsix$. For $m \in \ZZ \pm \frac{1}{4}$ and $m>1$, let 
$$
\sigma_e(m) = (q^{\delta}, q^{\delta+1}, \ldots, q^{2m-\frac{3}{2}}).
$$
Define the residual points $r_e(m)$ recursively by putting 
$$
r_e(\frac{1}{4}) = r_e(\frac{3}{4}) := \emptyset; \quad r_e(m) = (\sigma_e(m); r_e(m-1) ) \text{ if } m>1. 
$$
Finally we define the representing morphism $\xi_{T, \emm, \emp}: T^n \to T^l$ of the extra-special STM by 
\begin{equation}\label{eq:extramorphism}
\xi_{T, \emm, \emp}(t_1, \ldots, t_n) = \big(-r_e(\emm), v^{-1}t_1, vt_1, \ldots, v^{-1}t_n, vt_n, r_e(\emp) \big). 
\end{equation}

\end{enumerate}

\section{The proofs}\label{sec:pf}

We will prove in this section that the morphisms $\phi_{T, \emm, \emp}$ \eqref{eq:stmphi}, $\psi_{T, \emm, \emp}$ \eqref{eq:stmpsi} and $\xi_{T, \emm, \emp}$ \eqref{eq:extramorphism} satisfy the condition (T3) in the definition of spectral transfer morphisms, and hence they are indeed spectral transfer morphisms (the conditions (T1), (T2) are easy to verify, and we will not do it). A key idea to the proof is the factorisation of the $\mu$-function by a parabolic root sub-system (cf.~Section \ref{sec:IJKmu}). We shall interpret this parabolic factorisation in terms of Lusztig's classification of unipotent representations.

\subsection{Lusztig's classification and the parabolic factorisation of the $\mu$-function}\label{sec:IJKmu}

Recall that Lusztig's classification \cite{Lusztig} is achieved by matching the Iwahori--Matsumoto presentations of two affine Hecke algebras. One of these affine Hecke algebras is $\ha_{\utype} = \ha_{upt}(\PH^{F_u}, \sigma)$, derived from the irreducible unipotent representation $\pi \in \mathrm{Irr}(\mathbf{G}(k), \PH^{F_u}, \sigma)$ (as we have seen in \S 2.1). The other one, here we denote it by $\mathsf{H}:=\mathsf{H}(s, \mathcal{N, F})$, is obtained by considering cuspidal local systems $\mathcal{F}$ on unipotent orbits $\mathcal{N}$ in a connected complex reductive group $M:=Z_{G^\vee}(s)^\circ$ where $s \in {}^L G$ is semisimple and $G^\vee$ is the complex dual group of $G^{F_u}$. The conclusion is, $\pi$ is corresponding to the data $(s, \mathcal{N, F})$ if and only if the two affine Hecke algebras $\mathsf{H}(s, \mathcal{N, F})$ and $\ha_{upt}(\PH^{F_u}, \sigma)$ have the same Iwahori--Matsumoto presentation. 

Let $I$ denote a set of simple roots of $G^\vee$, and let $J\subset I$ be the subset corresponding to the reductive subgroup $M$. One can establish a bijection between $K=I-J$ and the simple roots in the underlying root system of $\mathsf{H}$. (The case that $K=\emptyset$ is very simple, and we shall not consider this case here.) More precisely, starting from the root systems of types $I$ and $J$, one can determine a unique reduced root system whose base is in bijection with $K$. This reduced root system is the underlying root system of $\mathsf{H}$. Then, the Iwahori--Matsumoto presentation of $\mathsf{H}$ is determined, with the explicit expression of the parameters given in \cite[5.12]{Lusztig}. The same procedure applies to $(\tisfI-\mathsf{J})/\theta_u$ and hence the Iwahori--Matsumoto presentation of $\ha_{upt}(\PH^{F_u}, \sigma)$ follows. (See \cite[1.17--1.20]{Lusztig}.) In particular, via Lusztig's classification, we obtain a bijection between $K$ and the set $(\tisfI-\mathsf{J})/\theta_u$ of $\theta_u$-orbits of $(\tisfI - \mathsf{J})$.


Now we turn to spectral transfer morphism. We consider the following morphism:
$$
\mathtt{Spec}\, Z(\ha_{upt}(\PH^{F_u}, \sigma)) \to \mathtt{Spec}\, Z(\mathsf{H}_{I}) =: T^{I},
$$ 
where $\mathsf{H}_{I}$ is an affine Hecke algebra attached to the root system $I$ (and given by the same recipe as the derivation of $\mathsf{H}$). By Lusztig's classification, $\ha_{upt}(\PH^{F_u}, \sigma)$ can be replaced by $\mathsf{H}$, so we can turn to consider a morphism between complex tori
$$
\mathcal{S}: T^K \to T^{I},
$$ 
where $T_K$ is the complex torus associated with $\mathsf{H}$. Let $T^J$ be a subtorus of $T^{I}$ determined by the subset $J\subset I$ (that is, if a root $\alpha \in \RR J \cap I$, then $\alpha | T^J$ is constant). We now interpret the image of $\mathcal{S}$ as a coset $L$ of $T^J$. In order that $\mathcal{S}$ to be a spectral transfer morphism, we must choose an appropriate base point $\respt \in T^J \cap (T^{I} \slash T^J)$ of $L$ to verify the condition (T3). To this end, we are led to the conclusion that the unitary part of $\respt$ and the semisimple element $s$ in the triple $(s, \mathcal{N, F})$ must be $\weylh_0$-conjugate (cf.~\cite[Appendix B]{Opd3}, \cite{Reexceptional,Opdl,FeOp} and also \cite[Conj.~1.5]{HII}).  

In this regard, the three morphisms \eqref{eq:stmphi}, \eqref{eq:stmpsi} and \eqref{eq:extramorphism} explicitly give the coordinates of the point $\respt$, as a residual point of the $\mu$-function attached to $\mathsf{H}_{I}$. Further, those coordinates in $\respt$ not involving with the original coordinates in $T_K$, give rise to a residual point of the root system of the type of $J$ (in the sense of \cite{OpdSol}). 

In general, suppose $L \subset T$ is a coset of a subtorus $T^L \subset T$. The natural projection $T \to T_L:=T/T^L$ corresponds to a sublattice $X_L \subset X$.  Then $L$ is a principal homogeneous space of $T^L$. Suppose $\respt$ is the base point of $L$, then we can write $L=\respt T^L$.

The intersection $R_L:=R_0 \cap X_L$ is a parabolic root sub-system of $R_0$. If necessary, we use an appropriate element in $\weylh_0$ to map $R_L$ to a standard position, then, we can choose a base $J$ of $R_L$ which is a subset of the base $F_0$ of $R_0$. We have a root datum $\calr_L:=(X^*(T_L), R_L, X_*(T_L), R_L^\vee)$ and a parabolic subgroup $\weylh_0^L$ of $\weylh_0$. (In fact, the root datum $\calr_L$ and the Weyl group $\weylh_0^L$ only depend on the subset $J$ of $F_0$.) Consequently, we can decompose the $\mu$-function as 
\begin{equation}\label{eq:parabolicmu}
\mu_{\zm, d}^{(L)} (\calr) = \mu_{\zm_L, d}^{\{\respt\}}(\calr_L) \frac{v^{-2\zm_{\weylh_0} (w^L)} }{\prod_{\alpha \in R_{0, +}\backslash R_{L, +}} \left( c_\zm (\alpha) c_\zm(\alpha^{w^L}) \right)},
\end{equation}
where $w^L:=w_0 w_L^{-1} \in \weylh^L \subset \weylh$ is the longest element.


\subsection{Conventions and the proofs}

Before the proof, we have some conventions. Given a nonzero rational number $x$ we denote by $\epsilon(x)=x/|x| \in \{\pm 1\}$ its sign. Recall that the normalisation factor $d=d(v)$, regarded as a rational function in $v$, satisfies $d(v)=d(v^{-1})$, and hence so does the $\mu$-function. Thus, we do not need to care about the minus signs. Moreover, we will neglect the powers of $q$, because we are only interested in the $q'$-part. Consequently, we can replace every expression $1-q^{-\mathfrak{a}}$ by $1-q^{\mathfrak{a}}$. Accordingly, we can rewrite the $\mu$-function as
\begin{equation}\label{normalizedmufunc}
\mu_{m_{\pm}, d}(R) = \prod_{\alpha \in R_+} \frac{d (1-\alpha(t)^2)^2}{(1+ v^{2\emm} \alpha(t)) (1+v^{-2\emm} \alpha(t)) (1-v^{2\emp} \alpha(t)) (1-v^{-2\emp} \alpha(t))}
\end{equation}
To simplify the notations, in the proofs below, we write $\Res(\mu, \respt)$ for $\mu^{(\{ \respt \}) }(\respt)$. 

The involutions 
\begin{equation*}
\eta_+: \; \emm \mapsto -\emp, \emm\mapsto \emm, \quad
\eta: \text{ interchanging $\emp$ and $\emm$}
\end{equation*} 
act on the $\mu$-functions, inducing isomorphisms of normalised affine Hecke algebras. The group \textbf{Iso}=$\langle \eta_+, \eta \rangle$ is isomorphic to the dihedral group of order eight, and it is called the spectral isomorphism group (cf.~\cite[\S 7]{Opds}). By virtue of \textbf{Iso}, it is enough to verify the condition (T3) $\emp \geq \emm \geq 0$. If $\emp=1/2$ then $\phi_{T, \emm,1/2}=\eta_+$. Hence we may and will assume that $\emp>1/2$ from now on.

\begin{proof}[The morphism $\phi_{T, \emm, \emp}^r: T^r \to T^{r+\emp-1/2}$ \eqref{eq:stmphi}.]
We first check that the image $L$ of $\phi_{T, \emm,\emp}^r$
in $T^{r+\emp-1/2}$ is a residual coset. 
The subset $R_L \subset R_0^{(1)}$ (the affine extension of $R_0$) of roots which are constant on the image is the set
$$
R_L=\{t_i^{\pm1}: r<i\leq r+\emp-1/2\} \cup \{t_i t_j^{-1}: r<i\neq j\leq r+\emp-1/2\},
$$
a root system of type $\bbb_{\emp-1/2}$.
Consider the factorisation $T^{r+\emp-1/2}=T^r \times T^{\emp-1/2}$.
By \cite[Thm.~3.7(ii)]{Opds} 
we need to show that the point $\respt_{\emm,\emp-1}$
in the torus $\{\mathrm{id}\} \times T^{\emp-1/2}$ defined by
$t_i(\respt_{\emm,\emp-1})=v^{2(i-r)-1}$ for $i=r+1,\dots,r+\emp-1/2$
is a residual point with respect to the $\mu$-function
$\mu_L:=\mu_{\emm,\emp-1}^{\emp-1/2}$ associated to the roots
of $R_L$. By definition of the residual
points we need to check that there are precisely $\emp-1/2$
factors in the denominator
of $\mu_{\emm,\emp-1}^{\emp-1/2}$ which are zero in this point.
Indeed, these are the factors corresponding to the roots
$t_{i+1}t_i^{-1}$ (for $i=r+1,\dots,r+\emp-3/2$) and the factor
corresponding to $t_{r+\emp-1/2}$.

It is a straightforward computation to determine the factor $\Res(\mu_{m_L,d}(\calr_L),\respt_{\emm,\emp-1})$, which is equal to $d \cdot\Res(\mu_{m_L,1}(\calr_L), \respt_{\emm,\emp-1})$ (in the notation of (\ref{eq:parabolicmu}), with $d=d_{\emm, \emp-1}$). We obtain
\begin{equation}\label{eq:cst}
\Res(\mu_{m_L,1}(\calr_L), \respt_{\emm,\emp-1})=-q^{(\emm-1/2)(\emp-1/2)}
\frac{(q-1)^{{\emp-1/2}}}{\prod_{k=\emm - \emp +1}^{\emm + \emp -1}(q^k+1)}
\end{equation}
In this computation it is not necessary to compute the exact
power of $q$, since \textit{a priori} we know 
that $\Res(\mu_{m_L, d}(\calr_L), \respt_{\emm,\emp-1})$ is invariant with respect to the substitution $v\mapsto v^{-1}$.
This remark simplifies the computation considerably.

In order to establish (T3) for $\phi_T$ we need to check firstly that
\begin{equation}
d_{\emm,\emp} =\lambda \cdot [d_{\emm,\emp-1} \cdot \Res(\mu_{m_L,1}(\calr_L), \respt_{\emm,\emp-1}) ]
\end{equation}
for some $\lambda\in\QQ^\times$. This is easy using (\ref{eq:cst}).
(Again, it is not necessary to compute the power of $q$).

Next, we need to establish (T3) for the non-constant part of $\mu^L$, by considering the restriction to $L$ of the roots in $\bbb_{r+\emp-1/2}$ which are not constant on $\{\mathrm{id}\} \times T^{\emp-1/2}$.
For the long roots this is trivial, and by the Weyl group symmetry
of $\mu$ we only need to consider the roots that restrict to
one fixed short root and its opposite, let us say
$t_{r+1}$ and $t_{r+1}^{-1}$. It is straightforward to collect all such
factors and check that, after cancellation of equal terms in numerator and
denominator, we will obtain
\begin{equation}
\frac{(1-t_{r+1}^{-2})(1-t_{r+1}^{2})}
{(1+v^{-2\emm}t_{r+1}^{-1})(1+v^{-2\emm}t_{r+1})(1-v^{-2\emp}t_{r+1}^{-1})(1-v^{-2\emp}t_{r+1})}.
\end{equation}
This is a factor of $\Res(\mu_{m_L, 1}(\calr_L), \respt_{\emm, \emp-1})$ with $m_L = (\emm, \emp)$. In short, the pull-back of the STM $\phi_{T, \emm,\emp}$ \eqref{eq:stmphi} on the $\mu$-function $\mu_{\emm, \emp-1, d_{\emm, \emp-1}}$, is equal to the $\mu$-function $\mu_{\emm, \emp, d_{\emm, \emp}}$ (up to a non-zero rational constant). As desired, finishing the proof. 
\end{proof}

Next, we discuss the case that $\zm_\pm \in \mathbb{Z}$ and $\emp\neq 0$.

\begin{proof}[The morphism $\psi_{T, \emm,\emp}^r: T^r \to T^{r+2\emp-2}$ \eqref{eq:stmpsi}.] We will sketch a proof similar to the proof of $\phi_{T, \emm, \emp}$. It is enough to consider the case $\emp>0$.
Let us check that the image of $\psi_{T, \emm, \emp}^{r}$ is a residual coset.
The parabolic subset $R_L$ of roots which are constant on the
image is the set
$$
R_L=\{t_i^{\pm1}: r<i\leq r+2(\emp-1)\}\cup\{t_it_j^{-1}: r<i\neq j\leq r+2({\emp}-1)\},
$$
a root system of type $\bbb_{2({\emp}-1)}$.
Consider the factorisation $T^{r+2({\emp}-1)}=T^r\times T^{2({\emp}-1)}$.
Again by \cite[Theorem 3.7(ii)]{Opds}
we need to show that the point $\respt_{\emm,\emp-2}$
in the torus $\{\text{id}\}\times T^{2({\emp}-1)}$ defined by
$t_i(\respt_{\emm,\emp -2})=v^{2\lfloor(i-r)/2\rfloor}$ for $i=r+1,\dots,r+2({\emp}-1)$
is a residual point with respect to the $\mu$-function
$\mu_L:=\mu_{\emm,\emp-1}^{r+2({\emp}-1)}$ associated to the roots
of $R_L$. Precisely $2({\emp}-1)$ roots are equal to $1$ on this point,
so by the definition of residual points we need to check that 
there are precisely
$2({\emp}-1)+2({\emp}-1)=4({\emp}-1)$ factors in the denominator
of $\mu_{\emm,\emp-1}^{2({\emp}-1)}$ which are zero in this point.
Indeed, there are $4({\emp}-2)$ such factors corresponding to
roots of the form $t_it_j^{-1}$, $2$ such factors from the
roots $t_i$ with $i=r+2{\emp}-4$ and $i=r+2{\emp}-3$, and finally
$2$ such factors from the roots $t_{r+1}t_{r+2}$ and
$t_{r+1}t_{r+3}$.

This time we have (in the notation of (\ref{eq:parabolicmu}),
with $d=d_{\emm,\emp-2}$):
\begin{equation}\label{eq:cstsympl}
\Res(\mu_{m_L,1}(\calr_L), \respt_{\emm,\emp-2})=
\frac{q^{\emm(\emp-1)}(q-1)^{2(\emp-1)}}
{\prod_{k=\emm-\emp+2}^{\emm+\emp-2}(q^k+1)^2(q^{(\emp-\emm)-1}+1)(q^{(\emp+\emm)-1}+1)}
\end{equation}
Again, in this computation it is not necessary to compute the exact
power of $q$.

In view of the factorisation (\ref{eq:parabolicmu}) of $\mu^L$,
in order for (T3) to hold we need to check firstly that
\begin{equation}
d_{\emm,\emp}=
\lambda \cdot [d_{\emm,\emp-2} \cdot 
\Res(\mu_{m_L,1}(\calr_L), \respt_{\emm,\emp-2})]
\end{equation}
for some $\lambda\in\QQ^\times$. This is easy using (\ref{eq:cst}).
(Again, it is not necessary to compute the exact power of $q$). The rest of the argument is an elementary
computation as in the proof for $\phi_{T, \emm, \emp}$, finishing the proof.
\end{proof}

Now we turn to the most complicated case. \\
\begin{prop}\label{verify-es-stm}
The morphism $\xi_{T, \emm, \emp}^{r}$ defines a spectral transfer morphism 
\begin{equation}\label{e-sstm}
\xi_{T, \emm, \emp}^{r}: T^r \to T^l
\end{equation}
with $l$ as defined in \eqref{eq:l}, which is called extra-special.\\
\end{prop}

\begin{rk}
Note that the unipotent affine Hecke algebra $\ha^r_{\emm, \emp}$ associated with $T^r$ has parameter $q^2$, while the unipotent affine Hecke algebra $\ha^{l}_{\delta_-, \delta_+}$ associated with $T^l$ has parameter $q$.
\end{rk}

We will use induction, and the following lemma will serve as the induction basis.\\
\begin{lemma}\label{lem:indbase}
Proposition \ref{verify-es-stm} is true for $(\emm,\emp)=(1/4,1/4),(3/4,1/4),(1/4,3/4)$ and $(3/4,3/4)$.
\end{lemma}

\begin{proof}[Proof of Lemma \ref{lem:indbase}]
Suppose $m_{\pm}=1/4$. Then $\delta_{\pm}=0$, and $a=b=0$, and $r_e(1/4) = \emptyset$. If moreover the rank $r=0$ then we have nothing to verify because $d_{1/4, 1/4}$ reduces to $1$, as does $d_{0,0}$, and there is no root contributing the $\mu$-functions so that both 
$\mu$-functions reduce to $1$ as well.

Now assume $r > 0$. Since $m_{\pm}=1/4$, we note that this is of parameter type VI with $a=b=0$. 
Hence $l=2r$, so that  
$d_{1/4, 1/4}=(v^2-v^{-2})^{-r}$ (recall that we have parameter $q^2$ here) and $d_{0,0}=(v-v^{-1})^{-2r}$. 
The $\mu$-function $\mu_{0,0}^{2r}$ has only contributions from the type $\ddd$ roots: 
 \begin{equation*}
\mu^{2r}_{0,0} = (v-v^{-1})^{-2r}\prod_{1 \leq i<j \leq 2r} \frac{(1-t_it_j^{-1})^2 (1-t_it_j)^2}{(1-qt_it_j^{-1}) (1-q^{-1}t_it_j^{-1}) (1-qt_it_j) 1-q^{-1}t_it_j)}.
\end{equation*}

We compute the pull-back $\xi^{\ast}(\mu^{2r}_{0,0})^{(L)}$, so we need to substitute 
$(t_1,t_2,\dots,t_{2r-1},t_{2r})$ by 
$$
(vs_1,v^{-1}s_1,vs_2,v^{-1}s_2,\dots,vs_{r},v^{-1}s_{r})
$$ 
in $\mu^{2r}_{0,0}$, 
after regularising the expression along the image $L$ of $\xi$. Concerning the regularisation, observe that the parabolic subsystem of roots which is constant on the image is 
$$
\{(t_1t_2^{-1})^{\pm 1},(t_3t_4^{-1})^{\pm 1},\dots,(t_{2r-1}t_{2r}^{-1})^{\pm 1} \},
$$
of type $\aaa_1^r$. After dropping the singular factors, the remaining ``constant factors" yield (including the normalisation factor $(v-v^{-1})^{-2r}$), up to a power of $v$ which is irrelevant for us (as explained before): $(v^2-v^{-2})^{-r}$, which is indeed the normalisation factor $d_{1/4, 1/4}$.  

For the non-constant part, by Weyl group invariance, it is enough to check the contribution of one type $\bbb$ root in each $W_0$-orbit.
First consider $s_1s_2$. The roots whose pull-back along $\xi$ equals a nonzero power of $s_1s_2$ times a power of $v$ are $t_1t_3, t_2,t_3, t_1t_4,t_2t_4$ (together with their opposites). These yield $qs_1s_2,s_1s_2,s_1s_2$ and $q^{-1}s_1s_2$ respectively, so these root give us the factor (after some cancellations, and up to powers of $v$ and other characters of $T^r$, which are irrelevant anyhow):
\begin{equation}
\frac{(1-s_1s_2)^2}{(1-q^{-2}s_1s_2)(1-q^2s_1s_2)}
\end{equation}
which shows that the pull back of $\mu_{0,0}^{2r}$ yields a factor 
which is a type $\ddd$ $\mu$ function with base $q^2$, as desired for $\mu^r_{1/4,1/4}$. 
 
The type $\ddd$ roots which pull back to a power of $v$ times a nonzero power of $s_1$ (a type $A_1$ root) 
are only $t_1t_2$ (and its opposite). Its pull back is $s_1^2$. 
This gives a factor in $\xi^{\ast}(\mu^{2r}_{0,0})^{(L)}$ of the form: 
\begin{equation}
\frac{(1-s_1)^2 (1+s_1)^2}{(1+v^{-1}s_1)(1+vs_1)(1-v^{-1}s_1) (1-vs_1)}
\end{equation}
This second fraction is exactly the factor in $\mu^r_{1/4,1/4}$ (with parameter $q^2$) 
for this last remaining type of root in the type $\bbb$ root system.
In other words, the condition (T3) is satisfied. Hence we have verified that $\xi_{T,1/4,1/4}^r$ 
represents an extra-special STM. 
\end{proof}

Now we turn to the proof of Proposition \ref{verify-es-stm} itself.

\begin{proof}[Proof of Proposition \ref{verify-es-stm}]

Recall that $m_{\pm} \in \ZZ \pm \frac{1}{4}$. Using spectral isomorphisms from $\text{Iso} = \langle \eta_+, \eta \rangle$, we can assume that $m_{\pm} >0$. 
So we can write $\emm = \kam + (2\epm -1)/4$ with $\kam \in \ZZ_{\geq 0}$ and $\epm \in \{0 ,1\}$. Let $\dem$ be defined by $\kam - \dem \in 2\ZZ$. 
In other words, $\dem =0$ (resp.~1) if $\kam$ is even (resp.~odd). Similarly we have $\emp = \kap + (2\epp -1)/4$ and $\dep$.

We will apply an inductive argument on $\emm+\emp$, where the induction base is provided by Lemma \ref{lem:indbase}. Let us first assume 
that $\emp\geq \emm>0$ and that $\emp>1$. By induction we may now assume that $\xi_{T,\emm,\emp-1}^r$ represents an STM.

Note that $\delta_{\pm}, \epsilon_{\pm} \in \{0, 1\}$. 
We define $\delta_{\pm}^c, \epsilon_{\pm}^c$ by the rules that 
$\delta_{\pm} + \delta_{\pm}^c =1, \epsilon_{\pm} + \epsilon_{\pm}^c =1$. 
Also we define $A(m)=2m-\frac{3}{2}$ and write 
$A_{\pm} = A(m_{\pm}) = 2(\kappa_{\pm}-1) + \epsilon_{\pm}\in\mathbb{Z}_{\geq -1}$. Observe that $A_+\geq 1$ by our assumptions.

To proceed, we need more notations. The $\mu$-function associated to the source normalised Hecke algebra $(\ha_{\emm,\emp}^r, \tau^d)$ 
of the alleged STM $\xi=\xi^r_{\emm,\emp}$ with parameters $m_{\pm}$ 
will be denoted by $\mu_{\emm, \emp , d}^r$, where we often omit the rank $r$ if there should be no confusion. 
When $d=d_{\emm,\emp}^r$ we  will simply write $\mu_{\emm, \emp}^r$. 
The $\mu$-function of the target is denoted by $\mu_{\delta_-, \delta_+}^l$, with $l$ given as  in (\ref{eq:l}).
Recall that we have, up to irrelevant factors, 
$\mu_{\delta_-, \delta_+}^l=d_{\delta_-,\delta_+}^l\mu^{\ddd,l}\mu_{\delta_-, \delta_+}^{\aaa,l}$ with 
\begin{equation}
d_{\delta_-,\delta_+}^l=(v-v^{-1})^{-l}(v+v^{-1})^{-\delta_-\delta_+}.
\end{equation}
The second factor arises in \cite[Prop.~2.5]{Opdl} from the fact that the reductive quotient $\overline{\mathbb{P}}$ of a minimal $F$-stable parahoric in the case $(\delta_-,\delta_+)=(1,1)$ is an $\mathbb{F}_q$-torus of split rank $l$ whose maximal $\mathbb{F}_q$-anisotropic subtorus has $q+1$ rational points over $\mathbb{F}_q$. Furthermore, 

\begin{equation}
\begin{aligned}
\mu^{\ddd,l} &= 
\prod_{1 \leq i < j \leq l} \frac{(1-t_it_j^{-1})^2}{(1-qt_i t_j^{-1})(1-q^{-1} t_i t_j^{-1})} \prod_{1 \leq i< j \leq l} \frac{(1-t_i t_j)^2}{(1-qt_i t_j) (1-q^{-1} t_i t_j)},\\
\mu_{\dem, \dep}^{\aaa,l} &= \prod_{i=1}^l \frac{(1-t_i^2)^2}{(1+q^{\dem} t_i) (1+q^{-\dem} t_i) (1-q^{\dep} t_i) (1-q^{-\dep} t_i)}.
\end{aligned}
\end{equation}

On the other hand $\mu_{\emm, \emp}^r=d_{\emm,\emp}^r\mu^{\ddd,r}(q^2)\mu_{\emm,\emp}^{\aaa,r}$ with 
$$
d_{\emm, \emp}^r=(v-v^{-1})^{-r}d^0_{\emm,\emp},
$$ 
and the normalisation factor $d^0_{\emm,\emp}$ is given by 
\begin{equation*}
d^0_{\emm,\emp}= 
\prod_{i=1}^{\lfloor \emm-\emp \rfloor} \big(\frac{v^{2|\emm - \emp| -2i}}{1+q^{2|\emm-\emp|-2i}} \big)^i 
\prod_{j=1}^{\lfloor \emm+\emp \rfloor} \big(\frac{v^{2|\emm + \emp| -2j}}{1+q^{2|\emm+\emp|-2j}} \big)^j
\end{equation*}
as in \cite[Eq.~(36)]{Opdl}. Here $\mu^{\ddd,r}(q^2)$ is similar to $\mu^{\ddd,l}$, only with rank $r$ instead of $l$ and with parameter $q^2$ instead of $q$. Finally:
\begin{equation}
\mu_{\emm, \emp}^{\aaa,r} = 
\prod_{i=1}^r \frac{(1-s_i^2)^2}{(1+v^{4\emm} s_i) (1+v^{-4\emm} s_i) (1-v^{4\emp} s_i) (1-v^{-4\emp} s_i)}
\end{equation}

Recall that the morphism $\xi:=\xi_{T,\emm,\emp}^r: T_r \to T_l$ is defined by
\begin{equation}\label{eq:xi}
\xi_{T,\emm,\emp}^r(s_1, \ldots, s_r) = (-r_e(\emm), v^{-1}s_1, vs_1, \ldots, v^{-1}s_r, vs_r, r_e(\emp) )
\end{equation}
where $r_e(\frac{1}{4}) = r_e(\frac{3}{4}) = \emptyset, r_e(m) = (\sigma_e(m), r_e(m-1) )$ and 
$$
\sigma_e(m) = (q^{\delta}, q^{\delta+1}, \ldots, q^{2m-\frac{3}{2}}).
$$ 
We denote $\respt_0 := (-r_e(\emm), r_e(\emp))$. 
\textbf{Note} that if we re-order the coordinates of $\xi(s_1, \ldots, s_r)$ 
(or of $\respt_0$), or invert them, then the result lies in the same $W_{2,0}$-orbit.
The $\mu$-function is invariant under the $W_{2,0}$-action, hence also invariant under such operation.\\

We now consider the rank 0 case, i.e.~$r=0$. This is the main challenge, as we will see. 
 
We need to verify the condition (T3), assuming that $\emp \geq \emm>0$ and 
$\emp>1$, and (by the induction hypothesis) that $\xi_{T,\emm,\emp-1}^0$ represents an STM.

Write $l_0=\lfloor l_-\rfloor+\lfloor l_+\rfloor$ for the rank of the target Hecke algebra if we use the 
parameters $(\emm, \emp)$, and $l_0'=\lfloor l_-\rfloor+\lfloor l_+'\rfloor$ if we use 
the parameters $(\emm',\emp'):=(\emm, \emp-1)$. Observe that $\dem'=\dem$ and $\dep'=\dep^c$.
By our normalisations of the Hecke algebras it suffices to show the following identity 
for the ratio of residues:  
\begin{equation}\label{eq:factorC}
C_{\emm, \emp -1} := 
\frac{\Res(\mu^{\ddd,l_0}, \respt_0)\cdot \Res(\mu_{\dem, \dep}^{\aaa, l_0}, \respt_0) }
{\Res(\mu^{\ddd, l_0'}, \respt_0' ) \cdot \Res(\mu_{\dem, \depc}^{\aaa, l_0'}, \respt_0')}
=\frac{d_{\dem,\depc}^{l_0'} \cdot \Res(\mu_{\dem, \dep}^{l_0}, \respt_0)}
{d_{\dem,\dep}^{l_0} \cdot \Res(\mu_{\dem, \depc}^{l_0'}, \respt'_0)}. 
\end{equation}
The last expression equals, up to powers of $v$ and rational constants,  
$$
\mathbf{A}(\emp) 
:= \frac{d_{\dem,\depc}^{l_0'} \cdot d^0_{\emm, \emp}}{d_{\dem,\dep}^{l_0} \cdot d^0_{\emm, \emp-1}}
=(v-v^{-1})^{A_++\delta^c_+}(v+v^{-1})^{\dem(\dep-\demc)} \frac{d^0_{\emm, \emp}} {d^0_{\emm, \emp-1}}.
$$
The second equality is easy to check, using the relations $l_0-l_0'=A_++\delta^c_+=2\emp-\frac{1}{2}-\delta_+$.

\textbf{Notation}. In the equations below we will simplify notations by omitting the references to the rank (if the arguments are given, the rank equals the number of coordinates of the argument so the explicit references to the ranks are superfluous), 
and we will simply write ``$\textup{reg}$" to indicate that we are using the regularisation of $\mu$-functions, i.e.~omitting the factors that are identically $0$ 
after evaluation at the argument. Finally, an expression like 
$\mu^{\ddd,\textup{reg}}(\respt_1;\respt_2)$
means that we only consider the product in the numerator and the denominator of those type $\ddd$-roots 
$t_i^{\pm 1}t_j^{\pm 1}$ for which $t_i$ is a coordinate of $\respt_1$ and $t_j$ is a coordinate of $\respt_2$, and only those factors which are not identically $0$.

Since $\respt_0=(\respt'_0,\sigma_e(\emp))$, we see that $C_{\emm, \emp-1}$ is equal to 
\begin{equation}
\begin{aligned}
\frac{\mu^{\aaa,\textup{reg}}_{\dem, \dep}(\respt'_0)}{\mu^{\aaa,\textup{reg}}_{\dem, \depc}(\respt'_0)} 
&\mu_{\dem, \dep}^{\aaa,\textup{reg}}(\sigma_e(\emp))\times\\
&\mu^{\ddd,\textup{reg}}(\sigma_e(\emp))
\mu^{\ddd,\textup{reg}}(-r_e(\emm); \sigma_e(\emp))\mu^{\ddd,\textup{reg}}(\sigma_e(\emp); r_e(\emp-1))
\end{aligned}
\end{equation}

It is easy to see that 
\begin{equation}\label{eq:deltachange}
\frac{\mu^{\aaa,\textup{reg}}_{\dem, \dep}(\respt'_0)}{\mu^{\aaa,\textup{reg}}_{\dem, \depc}(\respt'_0)} = \prod_{u_i \in r_e(\emm)} \frac{(1+q^{\depc} u_i) (1+q^{-\depc} u_i)}{ (1+q^{\dep} u_i) (1+q^{-\dep} u_i)} 
\prod_{w_j \in r_e(\emp-1)} \frac{(1-q^{\depc} w_j) (1-q^{-\depc} w_j)}{ (1-q^{\dep} w_j) (1- q^{-\dep} w_j)} 
\end{equation}
while 
\begin{equation}
\begin{aligned}
\mu^{\aaa,\textup{reg}}_{\dem, \dep} (\sigma_e(\emp) ) &= \prod_{t_k \in \sigma_e(\emp)} \frac{(1-t_k^2)^2}{(1+q^{\dem} t_k) (1+q^{-\dem} t_k) (1-q^{\dep} t_k) (1-q^{-\dep} t_k)}\\
& =\begin{cases}
1 & (\dem, \dep)=(0, 0)\\
\frac{1+q^{2\emp -3/2}}{(1+q)(1+q^{2\emp-1/2})} & (\dem, \dep)=(1, 0)\\
\frac{(1-q^{2\emp -3/2})(1-q)}{1-q^{2\emp-1/2}} & (\dem, \dep)=(0, 1)\\
\frac{(1-q^{4\emp -3})(1-q^2)}{1+q^{2\emp-1/2}} & (\dem, \dep)=(1, 1)
	\end{cases} \\
& = \big(\frac{1+q^{A_+}}{1+q^{A_+ +1}} \big)^{\dem} \big(\frac{1-q^{A_+}}{1-q^{A_+ +1}} \big)^{\dep} (1+q)^{ (-1)^{\depc}\dem} (1-q)^{\dep}
\end{aligned}
\end{equation}
We denote this last expression by (P1) = $\mu^{\aaa,\textup{reg}}_{\dem, \dep}(\sigma_e(\emp))$. 

Next we consider $\mu^{\ddd, \textup{reg}}(-r_e(\emm); \sigma_e(\emp))$ and  $\mu^{\ddd,\textup{reg}} (\sigma_e(\emp); r_e(\emp-1))$. \begin{equation}\label{eq:muDlc}
\begin{aligned}
   \mu^{\ddd, \textup{reg}}&(-r_e(\emm); \sigma_e(\emp))\\
:=&\prod_{t_i \in \sigma_e(\emp),\; t_j \in -r_e(\emm)} \frac{(1-t_i t_j^{-1})^2 (1-t_i t_j)^2}{(1-qt_i t_j^{-1}) (1-q^{-1}t_i t_j^{-1}) (1-qt_i t_j) (1-q^{-1}t_i t_j)} \\
=&\prod_{t_j \in r_e(\emm)} \frac{(1+q^{\dep} t_j^{-1}) (1+q^{A_+} t_j^{-1}) (1+q^{\dep} t_j) (1+q^{A_+} t_j)}{(1+q^{\dep-1} t_j^{-1})(1+q^{A_+ +1} t_j^{-1})(1+q^{\dep -1} t_j)(1+q^{A_+ +1} t_j)}
\end{aligned}
\end{equation}
where $A_+ =2\emp -3/2$. We can likewise obtain 
\begin{equation}\label{eq:muDcr}
\begin{aligned}
 \mu^{\ddd, \textup{reg}}&(\sigma_e(\emp), r_e(\emp -1) )  \\
=& \prod_{t_i \in \sigma_e(\emp),\; t_j \in r_e(\emp -1)} \frac{(1-t_i t_j^{-1})^2 (1-t_i t_j)^2}{(1-qt_i t_j^{-1}) (1-q^{-1}t_i t_j^{-1}) (1-qt_i t_j) (1-q^{-1}t_i t_j)} \\
=& \prod_{t_j \in r_e(\emp -1)} \frac{(1-q^{\dep} t_j^{-1}) (1-q^{A_+} t_j^{-1}) (1-q^{\dep} t_j) (1-q^{A_+} t_j)}{(1-q^{\dep-1} t_j^{-1})(1-q^{A_+ +1} t_j^{-1})(1-q^{\dep -1} t_j)(1-q^{A_+ +1} t_j)}
\end{aligned}
\end{equation}

Observe that up to some power of $v$, we can cancel some factors from (\ref{eq:deltachange}), (\ref{eq:muDlc}) and (\ref{eq:muDcr}) and obtain 
\begin{equation}\label{eq:BandC}
\begin{aligned}
&\frac{\mu_{\dem, \dep}(\respt'_0)}{\mu_{\dem, \depc}(\respt'_0)} \times \mu(\ddd, q, -r_e(\emm), \sigma_e(\emp) ) \times \mu(\ddd, q,  \sigma_e(\emp), r_e(\emp -1) )  \\
=& \prod_{t_j \in r_e(\emm)} \frac{(1+q^{A_+} t_j^{-1}) (1+q^{A_+} t_j)}{(1+q^{A_+ +1} t_j^{-1})(1+q^{A_+ + 1} t_j)} \prod_{t_j \in r_e(\emp -1)} \frac{(1-q^{A_+} t_j^{-1}) (1-q^{A_+} t_j)}{(1-q^{A_+ +1} t_j^{-1})(1-q^{A_+ + 1} t_j)}. 
\end{aligned}
\end{equation}
Denote 
$$
\text{(P2)} = \prod_{t_j \in r_e(\emm)} \frac{(1+q^{A_+} t_j^{-1}) (1+q^{A_+} t_j)}{(1+q^{A_+ +1} t_j^{-1})(1+q^{A_+ + 1} t_j)} 
$$ 
and 
$$
\text{(P3)} = \prod_{t_j \in r_e(\emp -1)} \frac{(1-q^{A_+} t_j^{-1}) (1-q^{A_+} t_j)}{(1-q^{A_+ +1} t_j^{-1})(1-q^{A_+ + 1} t_j)}.
$$ 
They can be simplified further. Observe that
$$
r_e(\emm) = \bigg(\sigma_e(\emm), r_e(\emm -1) \bigg) = \bigg(\sigma_e(\emm), \sigma_e(\emm -1), \ldots, \sigma_e(\frac{7-2\epm}{4}) \bigg)
$$
The number of $\sigma_e$'s in $r_e(\emm)$ is  $\kam$. Recall that for $g \in \ZZ_{>0}$ we defined:  
$$
\sigma_e(g+ \frac{2\epm-1}{4} ) = (q^{\bar{g}}, q^{\bar{g} +1}, \ldots, q^{2(g-1) + \epm} )
$$ where $\bar{g}=0$ if $g$ is even and $\bar{g}=1$ if $g$ is odd. 
Therefore we can write (P2) as 
\begin{equation}
\begin{aligned}
\text{(P2)} &= \prod_{g=2-\epm}^{\kam} 
\prod_{t_j \in \sigma_e(g+ \frac{2\epm-1}{4})} \frac{(1+q^{A_+} t_j^{-1}) (1+q^{A_+} t_j)}{(1+q^{A_+ +1} t_j^{-1})(1+q^{A_+ + 1} t_j)} \\
&= \prod_{g=2-\epm}^{\kam}\frac{(1+q^{A_+ - (2g-2+\epm)}) (1+ q^{A_+ + \bar{g}})}{(1+q^{A_++1-\bar{g}}) (1+q^{A_+ + (2g-1 + \epm)})}\\
&= \prod_{g=2-\epm}^{\kam}\frac{(1+q^{A_+ - (2g-2+\epm)})}{(1+q^{A_+ + (2g-1 + \epm)})} \times  \frac{(1+ q^{A_+ + \bar{g}})}{(1+q^{A_++1-\bar{g}})}.
\end{aligned}
\end{equation}
Notice that $\dem$ indicates the parity of $\kam$, so
\begin{equation*}
\prod_{g=2-\epm}^{\kam} \frac{(1+ q^{A_+ + \bar{g}})}{(1+q^{A_++1-\bar{g}})} = 
\begin{cases}
1 & \text{ if } \epm \neq \dem\\
\frac{1+q^{A_+}}{1+q^{A_+ +1}} & \text{ if } \epm = \dem =0\\
\frac{1+q^{A_++1}}{1+q^{A_+}} & \text{ if } \epm = \dem =1.
\end{cases}
\end{equation*}
So (P2) is equal to
\begin{equation}\label{eq:B}
\bigg(\frac{1+q^{A_+ + \epm}}{1+q^{A_+ + 1 - \epm}} \bigg)^{\epm \dem + \epmc \demc} \frac{(1+q^{A_+ -2 + \epm}) (1+q^{A_+ -4 + \epm}) \cdots (1+q^{A_+ -A_-})}{(1+q^{A_+ + 3 - \epm}) (1+q^{A_+ + 5 - \epm}) \cdots (1+q^{A_+ + A_- +1})},
\end{equation}
where $A_- = 2\emm -3/2$. (Recall that $\emp \geq \emm>0$ and $\emp>1$, 
hence $A_+\geq A_-\geq -1$ and $A_+\geq 1$. Observe that $\text{(P2)}=1$ if $0<\emm<1$.)

Similarly we can compute that (P3) is equal to 
\begin{equation}\label{eq:C}
\bigg(\frac{1-q^{A_+ + \epp}}{1-q^{A_+ + 1 - \epp}} \bigg)^{\epp \depc + \eppc \dep} \frac{(1-q^{2}) (1-q^{4}) \cdots (1-q^{A_+ -2 + \epp})}{(1-q^{2A_+  -1}) (1-q^{2A_+ -3}) \cdots (1-q^{A_+ +3-\epm})}.
\end{equation}

The final term to compute is 
$\mu^{\ddd,\textup{reg}}(\sigma_e(\emp))=\mu^{\ddd_0,\textup{reg}}(\sigma_e(\emp))\mu^{\ddd_{\not=0},\textup{reg}}(\sigma_e(\emp))$.
Here 
$\ddd_0$ denotes the type $\ddd$-roots whose coordinates sum up to zero (this is a maximal proper parabolic root 
subsystem, irreducible of type $\aaa$), and $\ddd_{\not=0}$ denotes the remaining roots of type $\ddd$ (this is not a root subsystem). 

Recall that $\sigma_e(\emp) = (q^{\dep}, q^{\dep +1}, \ldots, q^{2\emp -3/2})$. We can easily compute that
\begin{equation}\label{eq:AinD} 
\mu^{\ddd_0,\textup{reg}}(\sigma_e(\emp)) = \prod_{1 \leq i < j \leq A_+ -\dep} \frac{(1-q^{j-i})^2}{(1-q^{j-i+1}) (1-q^{j-i -1})} =\frac{(1-q)^{A_++1 -\dep}}{(1-q^{A_+ +1 -\dep})}
\end{equation}
by considering the multiplicities of the range of $j-i \in \{ 1, 2, \ldots, A_+ -\dep -1\}$. The same idea applies to computing
$\mu^{\ddd_{\not=0},\textup{reg}}(\sigma_e(\emp))$:
\begin{equation}\label{eq:DinD}
\scriptsize
\begin{aligned}
&\mu^{\ddd_{\not=0},\textup{reg}}(\sigma_e(\emp)) = \prod_{1 \leq i < j \leq A_+-\dep} \frac{(1-q^{i+j+2(\dep -1)})}{(1-q^{i+j+2\dep-1}) (1-q^{i+j+2\dep-3})} \\
&= \begin{cases}\displaystyle
\frac{(1-q^{1+2\dep}) (1-q^{3+2\dep}) \cdots (1-q^{A_+ +\dep -1}) (1-q^{A_+ +\dep +1}) \cdots (1-q^{2A_+ -1})}{(1-q^{2\dep})(1-q^{2+2\dep}) \cdots (1-q^{A_+ +\dep -2}) (1-q^{A_+ +\dep +2}) \cdots (1-q^{2A_+})} & \text{ if } \epp = \dep\\
\displaystyle\frac{(1-q^{1+2\dep}) (1-q^{3+2\dep}) \cdots (1-q^{A_+ +\dep}) (1-q^{A_+ +\dep }) \cdots (1-q^{2A_+ -1})}{(1-q^{2\dep})(1-q^{2+2\dep}) \cdots (1-q^{A_+ +\dep -1}) (1-q^{A_+ +\dep +1}) \cdots (1-q^{2A_+})}  & \text{ if } \epp \neq \dep.
\end{cases}
\end{aligned}
\end{equation}
Here, if $\dep=0$ then the denominator starts with $(1-q^2)$. 

We denote 
$\mu^{\ddd,\textup{reg}}(\sigma_e(\emp))$, $\mu^{\ddd_0,\textup{reg}}(\sigma_e(\emp))$, $
\mu^{\ddd_{\not=0},\textup{reg}}(\sigma_e(\emp))$
respectively by (P4), (P4a) and (P4d).

Now we multiply (P1), (P2), (P3) and (P4). The 4 parameters $\epp, \dep, \epm, \dem$ take values in $\{0, 1\}$ independently. So basically we need to consider 16 cases (The parity of $A_+$ is the same as $\epp$). But we spot a simplification when taking the product of (P3) and (P4d). 
We see that: 
\begin{equation*}
\small
\begin{aligned}
&\text{(P3)} \times \text{(P4d)} \\
=& \bigg(\frac{1-q^{A_+ + \epp}}{1-q^{A_+ + 1 - \epp}} \bigg)^{\epp \depc + \eppc \dep} \frac{(1-q^{2}) (1-q^{4}) \cdots (1-q^{A_+ -2 + \epp})}{(1-q^{2A_+  -1}) (1-q^{2A_+ -3}) \cdots (1-q^{A_+ +3-\epm})} \times \text{(P4d)}\\
=& \frac{(1-q^{1+2\dep}) (1-q^{3+2\dep}) \cdots (1-q^{A_+ +1 -\epp})}{(1-q^{A_+ +2 + \epp}) \cdots (1-q^{2A_+})} \\
=& \frac{(1-q)[1-q^2](1-q^3)[1-q^4] \cdots (1-q^{A_+ + 1 -\epp})}{[1-q^2][1-q^4][1-q^6][1-q^8] \cdots (1-q^{2A_+})} \frac{(1-q^{A_+ +1})^{\epp}}{(1-q)^{\dep}} \\
=&  \frac{(1-q)[1-q^2](1-q^3)[1-q^4] \cdots (1-q^{A_+ + 1 -\epp})}{[1-q^2][1-q^4][1-q^6][1-q^8] \cdots (1-q^{2(A_+ + 1 -\epp)})} \frac{(1-q^{2(A_+ + 1)})^{\eppc} (1-q^{A_+ +1})^{\epp}}{(1-q)^{\dep}} \\ 
=& \frac{(1+q^{A_+ +1})^{\eppc}}{(1+q)(1+q^2) \cdots (1+q^{A_+ +1 -\epp})} \frac{1-q^{A_+ +1}}{(1-q)^{\dep}} \\
=& \frac{1}{(1+q)(1+q^2) \cdots (1+q^{A_+})} \frac{1-q^{A_+ +1}}{(1-q)^{\dep}}, 
\end{aligned}
\end{equation*}
no matter if $\epp$ equals to $\dep$ or not. Here in the third and the fourth equations we insert in both the numerators and denominators the factors in square brackets to produce the factors $1+q^{\star}\, (\star = 1, \ldots, A_+)$ in the denominator.

We proceed to combine with (P4a) and (P1) to obtain 
\begin{equation*}
\begin{aligned}
&\frac{(1-q)^{A_++1 -\dep}}{(1-q^{A_+ +1 -\dep})} \times \bigg( \frac{1-q^{A_+}}{1-q^{A_+ + 1}}\bigg)^{\dep}  \times \bigg( \frac{1+q^{A_+}}{1+q^{A_+ + 1}}\bigg)^{\dem} (1+q)^{\dem (-1)^{\depc}} (1-q)^{\dep} \\
 \times & \frac{1}{(1+q)(1+q^2) \cdots (1+q^{A_+})} \frac{1-q^{A_+ +1}}{(1-q)^{\dep}} \\
 =& \frac{(1+q)^{\dem (-1)^{\depc}} (1-q)^{A_++1 -\dep} }{(1+q)(1+q^2) \cdots (1+q^{A_+})} \bigg( \frac{1+q^{A_+}}{1+q^{A_+ + 1}}\bigg)^{\dem} \\
 =& \text{(P1)} \times \text{(P3)} \times \text{(P4)},
\end{aligned}
\end{equation*}
no matter the value of $\dep$.

Finally we multiply with (P2). Note first that 
\begin{equation*}
\bigg(\frac{1+q^{A_+ + \epm}}{1+q^{A_+ + 1 - \epm}} \bigg)^{\epm \dem + \epmc \demc}  \bigg( \frac{1+q^{A_+}}{1+q^{A_+ + 1}}\bigg)^{\dem} = \frac{1+q^{A_+ + \epm}}{1+q^{A_+ +1}}
\end{equation*} 
and hence no matter the values of $\epm$ and $\dem$, the total product $C_{\emm, \emp-1}$ is equal to
\begin{equation}\label{eq:multABCD}
\scriptsize 
\begin{aligned}
 \frac{(1+q)^{\dem \cdot (-1)^{\depc}} (1-q)^{A_++1 -\dep}}{(1+q)\cdots (1+q^{A_+}) (1+q^{A_+ +1})} &\frac{(1+q^{A_+ + \epm}) (1+q^{A_+ -2 + \epm}) \cdots (1+q^{A_+ - A_- })}{(1+q^{A_+ + 3 -\epm}) (1+q^{A_+ + 5 -\epm}) \cdots (1+q^{A_+ + A_- +1})} \\
= (1-q)^{A_+ + 1 -\dep} (1+q)^{\dem \cdot (-1)^{\depc}}&\times \frac{1}{(1+q)(1+q^2)\cdots (1+q^{A_+ - A_- -1})} \times \\
 &\times\frac{1}{  (1+q^{A_+ -A_- +1}) (1+q^{A_+ -A_- +3}) \cdots \cdots (1+q^{A_+ + A_- +1})}.
\end{aligned}
\end{equation}
(The second and third factor of this product may be equal to $1$, if $A_+\leq A_-+1$ or if $A_-=-1$ respectively). 
In view of the expression of $d_{\emm, \emp}^0$, we now verify easily that $C_{\emm, \emp-1}= \mathbf{A}(\emp)$. Furthermore, we remark that if $\emm >\emp$ then we should compute $C_{\emm-1, \emp}$ which can be obtained by changing the subscripts $_+$ in (P1) to (P4) to $_-$ and similarly to obtain $\mathbf{A}(\emm)$. One can likewise verify that $C_{\emm-1, \emp} = \mathbf{A}(\emm)$. To sum up, we have verified by induction on $\emm+\emp$ that in the rank $0$ case of $\xi_T$ indeed represents a spectral transfer morphism. 

Now we consider the \emph{positive rank} cases. Again, assume that $\emp\geq \emm$.  
Still we need to verify the condition (T3) for the $\mu$-functions. 
We put $s=(s_1,\dots,s_r)$. Let $\xi(s):=\xi_{T,\emm,\emp}^r(s)$ and 
$\xi'(s):=\xi_{T,\emm,\emp-1}^r(s)$ be given as in (\ref{eq:xi}). Observe that 
$\xi(s)=(\xi'(s),\sigma_e(\emp))$. Let $L$ be the image 
of $\xi$, and let $L'$ be the image of $\xi'$. Observe that the parabolic root system 
$R_L$ of roots of $R_0=\bbb_l$ which restrict to constant functions on $L$ is isomorphic to $\aaa_1^r\times \bbb_{l_0}$, 
while the roots of $R_0'=\bbb_{l'}$ restricting to constant roots on $L'$ form a parabolic subsystem of type $\aaa_1^r\times \bbb_{l'_0}$. Using \eqref{eq:factorC}
it suffices to prove that the ratio
\begin{equation}\label{eq:positiverank}
 \frac{\mu^{l,(L)}_{\dem, \dep}(\xi(s))}{\mu^{l',(L')}_{\dem, \depc}(\xi'(s))}
=\frac{(v+v^{-1})^{-r} \Res(\mu_{\dem,\dep}^{l_0}, \respt_0) \cdot \mu_{\dem,\dep,R_0\backslash R_L}(\xi(s))}
{(v+v^{-1})^{-r} \Res(\mu_{\dem,\depc}^{l_0'}, \respt_0') \cdot \mu_{\dem,\depc,R_0'\backslash R_{L'}}(\xi'(s))}
\end{equation}
which we denoted as $C^r_{\emm, \emp-1}(s)$, is equal to $C_{\emm, \emp-1} \times \frac{\mu_{\emm, \emp}^r(s)}{\mu_{\emm, \emp-1}^r(s)}$
where
$C_{\emm, \emp-1}$ is defined in \eqref{eq:factorC}, and
\begin{equation}\label{eq:positiverankmum}
\frac{\mu_{\emm, \emp}^r(s)}{\mu_{\emm, \emp-1}^r(s)} = \prod_{i=1}^{r} \frac{(1-q^{-2(\emp-1)} s_i) (1-q^{2(\emp-1)} s_i)}{(1-q^{-2\emp} s_i) (1-q^{2\emp} s_i)}.
\end{equation}
Based on the result of rank 0 case, to prove that $C^r_{\emm, \emp-1}(s)=C_{\emm, \emp-1} \times \frac{\mu_{\emm, \emp}^r(s)}{\mu_{\emm, \emp-1}^r(s)}$, we just need to verify that 
\begin{equation}\label{eq:verify}
\frac{\mu_{\dem,\dep,R_0\backslash R_L}(\xi'(s),\sigma_e(\emp))}{\mu_{\dem,\depc,R_0'\backslash R_{L'}}(\xi'(s))}=
\prod_{i=1}^{r} \frac{(1-q^{-2(\emp-1)} s_i) (1-q^{2(\emp-1)} s_i)}{(1-q^{-2\emp} s_i) (1-q^{2\emp} s_i)}.
\end{equation}
It is enough to consider the contribution to the factors of the right hand side involving $s_1$.
For this we need to consider the contribution in the numerator of the left hand side of the set of type $\ddd$-roots $t_i^{\pm 1}t_j^{\pm 1}$ such that 
$t_i=v^{-1}s_1$ or $vs_1$, and $t_j$ is a coordinate of $\sigma_e(\emp)$.
In addition, we need to consider the contribution in the numerator and 
denominator of the left hand side of the roots $t_1=v^{-1}s_1$, $t_2=vs_1$ and their opposites. Therefore, in this computation we may assume that $r=1$.

The contribution from the type $\aaa_1^r$ roots $t_1=v^{-1}s_1$ and $t_2=vs_1$ is:
\begin{equation*}
\begin{aligned}
 &\frac{\mu_{\dem, \dep}^2(v^{-1}s_1, vs_1)}{\mu_{\dem, \depc}^2(v^{-1}s_1, vs_1)} \\
 = &\frac{(1-q^{\depc -1/2} s_1) (1-q^{\depc +1/2} s_1) (1-q^{-\depc -1/2} s_1) (1-q^{-\depc +1/2} s_1)}{ (1-q^{\dep -1/2} s_1) (1-q^{\dep +1/2} s_1) (1-q^{-\dep -1/2} s_1) (1-q^{-\dep +1/2} s_1)} \\
 =& \bigg(\frac{(1-q^{3/2} s_1) (1-q^{-3/2} s_1)}{(1-q^{1/2} s_1) (1-q^{-1/2} s_1)} \bigg)^{(-1)^{\dep}}.
\end{aligned}
\end{equation*}
Next, consider the type $\ddd$-roots $t_1t_j^{\pm 1}$ and $t_2t_j^{\pm 1}$ in the numerator,  
with $t_1 = v^{-1}s_1$, $t_2=vs_1$, and $t_j \in \sigma_e(\emp)$. These yield:
\begin{equation*}
\begin{aligned}
& \prod_{s =\dep}^{2\emp -3/2} \frac{(1-q^{-s+1/2} s_1) (1-q^{-s-1/2} s_1) (1-q^{s+1/2} s_1) (1-q^{s-1/2} s_1)}{(1-q^{-s+3/2} s_1) (1-q^{-s-3/2} s_1) (1-q^{s+ 3/2} s_1) (1-q^{s-3/2} s_1)} \\
=& \frac{(1-q^{-2\emp+2} s_1) (1-q^{2\emp-2} s_1) (1-q^{-\dep -1/2} s_1) (1-q^{-\dep +1/2} s_1) }{(1-q^{-\dep +3/2} s_1) (1-q^{\dep -3/2} s_1) (1-q^{-2\emp} s_1) (1-q^{2\emp} s_1)}.
\end{aligned}
\end{equation*}
Multiplying the contributions from the type $\aaa$ and type $\ddd$ roots, we are quickly led to \eqref{eq:verify}. Using induction on $\emm+\emp$ (as in the rank zero case), 
and the induction basis Lemma \ref{lem:indbase}, we can finish the proof that $\xi^r_{T, \emm,\emp}$ represents an extra-special spectral transfer morphism.
\end{proof}

We conclude this paper by the following remark.
\begin{rk}
Compositions of $\phi_{T, \emm, \emp}$ and $\psi_{T, \emm, \emp}$ are obviously defined. By applying (the compositions of) $\phi_T, \psi_T$ and $\xi_T$, we obtain a spectral transfer morphism from a rank zero unipotent affine Hecke algebras to the one of maximal rank with parameters attain the minimum in the same parameter type. These three spectral transfer morphisms are necessary to prove the main result, Theorem 3.4 in \cite{Opdl}, called the \emph{essential uniqueness} of spectral transfer morphism. 
\end{rk}

\textbf{Acknowledgement} This work was supported by European Research Council (ERC) Advanced Grant [number 268105] and by Netherlands Organisation for Scientific Research (NWO) Vidi Grant [number 639.032.528]. The author wants to express his gratitude to Eric Opdam and Maarten Solleveld, for their patient guidance and many helpful discussions.

\end{document}